\documentclass[a4paper]{article}

\usepackage[dvipsnames]{xcolor}
\usepackage{amsmath, amsthm, amsfonts, bbold, float, tikz, ulem, mathtools, tabularx, mathrsfs, booktabs}
\usetikzlibrary{shapes, decorations.pathreplacing}

\RequirePackage{algorithm}[0.1] 
\usepackage{algpseudocode} 

\newcommand{\bm}[1]{\text{\boldmath $#1$\unboldmath}}
\newcommand{\bx}{\bm{x}}
\newcommand{\beff}{\bm{f}}
\newcommand{\bn}{\bm{n}}

\newcommand{\bu}{\bm{u}}

\newcommand{\bv}{\bm{v}}
\newcommand{\bU}{\bm{U}}
\newcommand{\bI}{\bm{I}}
\newcommand{\bb}{\bm{b}}

\newcommand{\bg}{\bm{g}}

\newcommand{\bV}{\bm{V}}

\newcommand{\bY}{\bm{Y}}

\newcommand{\bXi}{\bm{\Xi}}

\newcommand{\blambda}{\bm{\lambda}}

\newcommand{\bvarphi}{\bm{\varphi}}

\newcommand{\bmu}{\bm{\mu}}

\newcommand{\bLambda}{\bm{\Lambda}}

\newcommand{\Nmu}{N^{\mu}}

\newcommand{\pgd}{\texttt{PGD}}

\newcommand{\bsigma}{\bm{\sigma}}

\newcommand{\mat}[1]{\mathbf{#1}}

\newcommand{\N}{{\mathbb{N}}}  
\newcommand{\R}{{\mathbb{R}}}  

\newtheorem{tm}{Theorem}
\theoremstyle{remark}
\newtheorem{rem}[tm]{Remark}

\newcommand{\be}[1]{\textcolor{red}{#1}}

\newcommand{\all}{\text{\textit{:}}}

\graphicspath{{./figures/}}

\usepackage{a4wide}

\begin{document}

\begin{center}
\begin{Large}
\textbf{An overlapping domain decomposition method for parametric Stokes and Stokes-Darcy problems via proper generalized decomposition}
\end{Large}

\bigskip

Marco Discacciati$^1$, Ben J. Evans$^1$, Matteo Giacomini$^{2,3}$

\medskip

\begin{small}
${}^1$ Department of Mathematical Sciences, Loughborough University, Epinal Way, LE11~3TU, Loughborough, United Kingdom. m.discacciati@lboro.ac.uk, b.j.evans@lboro.ac.uk

${}^2$ Laboratori de C\`alcul Numeric (LaC\`aN), E.T.S. de Ingenier\'ia de Caminos, Canales y Puertos, Universitat Polit\`ecnica de Catalunya, Barcelona, Spain.
               
${}^3$ Centre Internacional de M\`etodes Num\`erics en Enginyeria (CIMNE), Barcelona, Spain. matteo.giacomini@upc.edu
\end{small}

\end{center}

\medskip

\begin{abstract}
A strategy to construct physics-based local surrogate models for parametric Stokes flows and coupled Stokes-Darcy systems is presented.
The methodology relies on the proper generalized decomposition (PGD) method to reduce the dimensionality of the parametric flow fields and on an overlapping domain decomposition (DD) paradigm to reduce the number of globally coupled degrees of freedom in space.
The DD-PGD approach provides a non-intrusive framework in which end-users only need access to the matrices arising from the (finite element) discretization of the full-order problems in the subdomains. 
The traces of the finite element functions used for the discretization within the subdomains are employed to impose arbitrary Dirichlet boundary conditions at the interface, without introducing auxiliary basis functions.
The methodology is seamless to the choice of the discretization schemes in space, being compatible with both LBB-compliant finite element pairs and stabilized formulations,
and the DD-PGD paradigm is transparent to the employed overlapping DD approach.
The local surrogate models are glued together in the online phase by solving a parametric interface system to impose continuity of the subdomain solutions at the interfaces, without introducing Lagrange multipliers to enforce the continuity in the entire overlap and without solving any additional physical problem in the reduced space.
Numerical results are presented for parametric single-physics (Stokes-Stokes) and multi-physics (Stokes-Darcy) systems, showcasing the accuracy, robustness, and computational efficiency of DD-PGD, and its capability to outperform DD methods based on high-fidelity finite element solvers in terms of computing times.
\end{abstract}



\section{Introduction}

The simulation of highly viscous flows and their interaction with permeable environments is a relevant problem in many fields, from environmental sciences (e.g., subsurface pollution transport)~\cite{Yotov-AKNY-19} and microfluidics (e.g., lab-on-chip design)~\cite{Doblare-OPRBASLMDOMO-24} to chemical engineering (e.g., filtration systems)~\cite{Discacciati-PBDD-19} and biomedical engineering (e.g., blood flow perfusion in tissues)~\cite{Zunino-DZ-11}.
Physical models describing these interactions of free-flow regions and porous media classically rely on coupled Stokes-Darcy models.

Whilst the computation of one instance of such problems has been extensively studied in the literature,  the case in which physical and material parameters, boundary conditions,  geometric configurations, and operating scenarios depend on user-defined parameters still represents a computational challenge for existing methodologies.
Difficulties of the resulting parametric Stokes-Darcy model include: (i) high-dimensionality of the parametric problems; (ii) large number of globally coupled degrees of freedom in space; (iii) multi-physics nature of the system.
It is worth noticing that effective solutions have been proposed to tackle (i) and (ii)+(iii) independently. The former rely on reduced order models (ROMs), e.g., reduced basis, proper orthogonal decomposition and, more recently, scientific machine learning to circumvent the \textit{curse of dimensionality}~\cite{Martini:2015:ACM,Benner-BF-15,Schillinger-SJNS-22,Hesthaven-DH-23}. The latter employ either non-overlapping or overlapping domain decomposition (DD) strategies to reduce the number of coupled unknowns in the physical domain~\cite{Discacciati-DQV-07,Gunzburger-CGHW-11,Yotov-VWY-14,Caiazzo-CJW-14,Discacciati-DG-18}.

This work aims to go beyond this dichotomy by concurrently exploiting the dimensionality reduction capabilities provided by proper generalized decomposition (PGD) to represent parametric flow fields~\cite{Chinesta:2014,Diez:2017:CMAME,RS-SBGH:20}, and the reduction of the size of the spatial problems obtained by means of  the interface control domain decomposition (ICDD) paradigm~\cite{Discacciati-DGGQ-16,Discacciati:2024:JCP}.
The proposed methodology builds upon the DD-PGD framework introduced in~\cite{Discacciati:2024:CMAME,Discacciati:2024:DD28} for parametric scalar convection-diffusion-reaction equations,  inheriting its appealing properties, while extending and generalizing multiple aspects of this approach. 
For the sake of completeness, note that DD for local surrogate modelling has been extensively studied in recent years in conjunction with \textit{a posteriori} ROMs, relying on both overlapping~\cite{Iollo-IST-23,Tezaur-MWGT-24} and non-overlapping~\cite{Rozza-PNTBR-23,Taddei-TXZ-24,Kuberry-HKB-25} strategies. Interested readers can also refer to~\cite{Buhr:2021,Klawonn-HKLW-21} for some recent reviews.

The DD-PGD strategy constructs a physics-based \textit{a priori} ROM as a rank-one approximation in each subdomain by means of a greedy approach that circumvents the need for sophisticated sampling strategies, see the discussion in~\cite{Giacomini-GBSH-21}.
Overlapping DD strategies are employed to decouple the degrees of freedom of the different subdomains, avoiding the introduction of a global parametric solution at the interface and the imposition of continuity of  fluxes required in non-overlapping approaches.
This allows the conception of a non-intrusive paradigm for the construction of local surrogate models, where the end-user only needs access to the matrices arising from the spatial discretization of the full-order problem in each subdomain, with arbitrary Dirichlet boundary conditions imposed at the interface.
The linearity of the operator is exploited to reduce the dimensionality of the local parametric problems via superposition.
Once the PGD local surrogate models are computed, they are evaluated in the online phase by simply performing interpolation in the parametric space, without the need to solve any additional problem in the reduced space, exploiting the non-intrusiveness of PGD-ROMs~\cite{Ladeveze-CNLB-16,Zou-ZCDA-18,Tsiolakis-TGSOH-20,Cavaliere-CZSLD-21,Tsiolakis-TGSOH-22,Cavaliere-CZSLD-22a,Cavaliere-CZSLD-22b}.
Finally, the local surrogate models are glued together solving a small interface problem in real time, without any call to the full-order solver.

In addition, the present work goes beyond the state of the art by presenting, for the first time,  PGD-based local surrogate models for a coupled, multi-physics problem, involving both scalar (i.e., pressure) and vector-valued unknowns (i.e., velocity).
More precisely, the DD-PGD method is extended to parametric Stokes problems and coupled Stokes-Darcy systems.
The methodology is completely transparent to the choice of the full-order solver and numerical experiments show that both LBB-compliant and stabilized finite element formulations, see~\cite{Quarteroni:1994,Donea:2003}, are suitable for the construction of accurate and stable local surrogate models.
Moreover,  the overlapping approach at the basis of the DD-PGD method only requires imposing equality of the local PGD velocity (or pressure) at the interface. This enables our approach to outperform existing PGD strategies based on domain decomposition by circumventing the need to introduce separated representations of either the Lagrange multiplier used to impose continuity of the solution on the entire overlap as in~\cite{Nazeer:2014:CM}, or the trace of the unknown solution as employed in~\cite{Huerta:2017:IJNME}.

The rest of this article is organized according to the following structure.
Section \ref{sect:Stokes} describes the multi-domain formulation of the parametric Stokes problem, extending the DD-PGD framework to the offline construction of local ROMs for incompressible flows and their online coupling via a surrogate-based alternating Schwarz method.
In Section \ref{sect:Darcy}, the local surrogate model for the mixed formulation of the parametric Darcy equation is introduced, together with a surrogate version of the ICDD algorithm for the online phase of the parametric Stokes-Darcy simulation.
Numerical experiments, reported in Section \ref{sect:numericalResults},  are presented to verify the accuracy and efficiency of the DD-PGD method for both parametric Stokes-Stokes and Stokes-Darcy problems. Moreover a cross-flow simulation is presented to showcase the capability of the proposed methodology to solve parametric studies of interest in the design of filtration systems.
Section \ref{sect:Conclusion} summarizes the results of this work and \ref{app:Details} reports some implementation details for the construction of the local surrogate models using the Encapsulated PGD Algebraic Toolbox~\cite{Diez:2020:ACME}.

\section{Parametric Stokes problem}
\label{sect:Stokes}

\subsection{Two-domain formulation of the parametric Stokes problem}
\label{sec:StokesProbStatement}

Consider an open bounded domain $\Omega \subset \R^d\ (d = 2, 3)$ with Lipschitz boundary $\partial\Omega$ partitioned into two non-empty disjoint parts $\Gamma^D$ and $\Gamma^N$, and with $\bn$ being the unit normal vector to $\partial\Omega$ pointing outwards of $\Omega$. Let $\bmu = (\mu_1, \dots, \mu_{N_p}) \in \mathcal{P}$ be a tuple of $N_p \in \N$ problem parameters with $\mathcal{P} = \mathcal{I}_1 \times \dots \times \mathcal{I}_{N_p} \subset \R^{N_p}$ and each  $\mathcal{I}_j$ compact ($j = 1, \dots, N_p$). The parametric Stokes problem describing the motion of a steady, viscous, incompressible fluid confined to $\Omega$ can be written as: for all $\bmu \in \mathcal{P}$, find the fluid velocity $\bu(\bmu)$ and pressure $p(\bmu)$ such that
\begin{subequations}\label{eq:StokesGlobal}
\begin{eqnarray}
		-\nabla \cdot \bsigma(\bu(\bmu), p(\bmu); \bmu) = \beff(\bmu) &&\quad \text{in } \Omega,\\
		\nabla \cdot \bu(\bmu) = 0 &&\quad \text{in } \Omega,\\
		\bu(\bmu) = \bg^D(\bmu) &&\quad \text{on } \Gamma^D,\\
		\bsigma(\bu(\bmu), p(\bmu); \bmu) \bn = \bg^N(\bmu) &&\quad \text{on } \Gamma^N,
\end{eqnarray}
\end{subequations}
where $\bsigma(\bu(\bmu), p(\bmu); \bmu) = 2\nu(\bmu)\nabla^s\bu(\bmu) - p(\bmu)\bI$ represents the Cauchy stress tensor with $\nabla^s\bu(\bmu) = \frac{1}{2}\left(\nabla\bu(\bmu) + (\nabla\bu(\bmu))^T\right)$. For all values of the problem parameters $\bmu \in \mathcal{P}$, $\beff(\bmu)$, $\bg^D(\bmu)$ and $\bg^N(\bmu)$ are regular-enough given functions, and the fluid viscosity is always positive, i.e., $\nu(\bmu) > 0$ in $\Omega$.

We decompose $\Omega$ into two overlapping subdomains $\Omega_1$ and $\Omega_2$: $\Omega = \Omega_1 \cup \Omega_2$, $\Omega_{12}:= \Omega_1 \cap \Omega_2 \neq \emptyset$, and, for $i=1,2$, let $\Gamma_i^D = \Gamma^D \cap \partial\Omega_i$ and $\Gamma_i^N = \Gamma^N \cap \partial\Omega_i$. To reformulate problem \eqref{eq:StokesGlobal} into an equivalent two-domain form, we assume that $\Gamma^N \cap \partial\Omega_{12} \not= \emptyset$ to fix the pressure in a unique way (see \cite{Discacciati:2013:Stokes} for details). Moreover, let $\Gamma_i = \partial \Omega_i \setminus \partial \Omega$ be the interfaces ($i=1,2$) that we assume do not intersect, i.e., $dist(\Gamma_1,\Gamma_2)>0$.

Then, a straightforward adaptation of the proof of Proposition 2.1 in \cite{Discacciati:2013:Stokes}, to account for the presence of the problem parameters $\bmu$, shows that the Stokes problem \eqref{eq:StokesGlobal} can be rewritten in the equivalent two-domain formulation: for all $\bmu \in \mathcal{P}$ and $i = 1, 2$, find the local fluid velocity $\bu_i(\bmu)$ and pressure $p_i(\bmu)$ such that
\begin{subequations}\label{eq:StokesMultiDomain}
	\begin{eqnarray}
		-\nabla \cdot \bsigma(\bu_i(\bmu), p_i(\bmu); \bmu) = \beff_i(\bmu) &&\quad \text{in } \Omega_i,\\
		\nabla \cdot \bu_i(\bmu) = 0 &&\quad \text{in } \Omega_i,\\
		\bu_i(\bmu) = \bg_i^D(\bmu) &&\quad \text{on } \Gamma^D_i,\\
		\bsigma(\bu_i(\bmu), p_i(\bmu); \bmu) \bn = \bg_i^N(\bmu) &&\quad \text{on } \Gamma^N_i,\\
		\bu_1(\bmu) = \bu_2(\bmu) &&\quad \text{on } \Gamma_1 \cup \Gamma_2, \label{eq:StokesMultiDomain_5}
	\end{eqnarray}
\end{subequations}
where $\beff_i(\bmu) = \beff(\bmu)_{\vert\Omega_i}$, $\bg_i^D(\bmu) = \bg^D(\bmu)_{\vert\Gamma_i^D}$ and $\bg_i^N(\bmu) = \bg^N(\bmu)_{\vert\Gamma_i^N}$.
The equivalence between \eqref{eq:StokesGlobal} and \eqref{eq:StokesMultiDomain} implies that, for all $\bmu \in \mathcal{P}$, $\bu_1(\bmu) = \bu_2(\bmu)$ and $p_1(\bmu) = p_2(\bmu)$ in $\Omega_{12}$, and that
\begin{equation}\label{eq:globalVelocityPressure}
 \bu(\bmu) =
 \left\{
 \begin{array}{ll}
 \bu_1(\bmu) & \text{in } \Omega_1 \\
 \bu_2(\bmu) & \text{in } \Omega_2 \setminus \Omega_{12}
 \end{array}
 \right.
 \quad \text{and} \quad
 p(\bmu) =
 \left\{
 \begin{array}{ll}
 p_1(\bmu) & \text{in } \Omega_1 \\
 p_2(\bmu) & \text{in } \Omega_2 \setminus \Omega_{12} \, .
 \end{array}
 \right.
\end{equation}

We consider now two regular computational meshes $\mathcal{T}_i$ in $\Omega_i$ ($i=1,2$) made by simplices or quadrilaterals/hexahedra. For simplicity of presentation, we assume that the meshes coincide in the overlapping region $\Omega_{12}$, and that they are conforming with the interfaces $\Gamma_1$ and $\Gamma_2$. For a suitable integer $r \geq 1$, let $\mathbb{P}_r$ and $\mathbb{Q}_r$ be the spaces of polynomials of degree less than or equal to $r$ with respect to each variable for simplices and quadrilaterals/hexahedra, respectively. 
Let $X_{i}^r = \{ v \in C^0(\bar{\Omega}_i) \, : \, v_{\vert K} \in \mathbb{P}_r \text{ or } v_{\vert K} \in \mathbb{Q}_r \text{ for all } K \in \mathcal{T}_i\}$, and, for $i=1,2$, define the finite dimensional spaces
\begin{equation*}
 \bV_i^h = \{ \bv \in [H^1(\Omega_i) \cap X_{i}^r]^d \; : \; \bv = \mathbf{0} \; \text{on } \Gamma_i^D \cup \Gamma_i \}, \quad W_i^h = L^2(\Omega_i) \cap X_{i}^s,
\end{equation*}
for suitable polynomial degrees $r,s \in \mathbb{N}$. Moreover, consider the discrete trace space $\bY_i^h = [Y_i^h]^d$ with
\begin{equation*}
  Y_i^h = \{ \lambda \in C^0(\Gamma_i) \, : \, \lambda = 0 \text{ on } \overline{\Gamma}_i^D \cap \overline{\Gamma}_i \text{ if } \overline{\Gamma}_i^D \cap \overline{\Gamma}_i \not= \emptyset, \text{ and } \exists \, v \in H^1(\Omega_i) \cap X_i^r \text{ s.t. } v = \lambda \text{ on } \Gamma_i \}.
\end{equation*}

Then, at the discrete level (but still writing the problems in strong form for the sake of clarity) in each subdomain $\Omega_i$, we introduce the following splitting of the local Stokes problem.
\begin{enumerate}
    \item Local Stokes problem depending on assigned data: find $\bu_{i,h}^f(\bmu)$ and $p_{i,h}^f(\bmu)$ such that
\begin{subequations}\label{eq:StokesMultiDomainF}
	\begin{eqnarray}
		-\nabla \cdot \bsigma(\bu_{i,h}^f(\bmu), p_{i,h}^f(\bmu); \bmu) = \beff_i(\bmu) &&\quad \text{in } \Omega_i,\\
		\nabla \cdot \bu_{i,h}^f(\bmu) = 0 &&\quad \text{in } \Omega_i,\\
		\bu_{i,h}^f(\bmu) = \bg_i^D(\bmu) &&\quad \text{on } \Gamma^D_i,\\
		\bsigma(\bu_{i,h}^f(\bmu), p_{i,h}^f(\bmu); \bmu) \bn = \bg_i^N(\bmu) &&\quad \text{on } \Gamma^N_i,\\
		\bu_{i,h}^f(\bmu) = \widetilde{\bg}_i^D(\bmu) &&\quad \text{on } \Gamma_i\,. \label{eq:StokesMultiDomainF_5}
	\end{eqnarray}
\end{subequations}
To avoid incompatible boundary conditions at the intersection of $\Gamma_i^D$ and $\Gamma_i$, the function $\widetilde{\bg}_i^D(\bmu)$ is introduced in \eqref{eq:StokesMultiDomainF_5}. More precisely, $\widetilde{\bg}_i^D(\bmu)$ is defined as 
\begin{equation}
    \widetilde{\bg}_i^D(\bmu) = \begin{cases}
        \widehat{\bg}_i^D(\bmu)\quad & \text{if } \overline{\Gamma}_i^D \cap \overline{\Gamma}_i \neq \emptyset,\\
        \mathbf{0} \quad &\text{if } \overline{\Gamma}_i^D \cap \overline{\Gamma}_i = \emptyset,
    \end{cases}
\end{equation}
where $\widehat{\bg}_i^D(\bmu)$ is a suitable continuous prolongation of $\bg_i^D(\bmu)$ onto $\Gamma_i$.
Moreover, let $\bg_{\Omega_i}^D(\bmu) \in [H^1(\Omega_i) \cap X_i^r]^d$ be a suitable continuous extension of the Dirichlet data $\bg_i^D(\bmu)$ and $\widetilde{\bg}_i^D(\bmu)$ such that $\bg_{\Omega_i}^D(\bmu) = \bg_i^D(\bmu)$ on $\Gamma^D_i$ and $\bg_{\Omega_i}^D(\bmu) = \widetilde{\bg}_i^D(\bmu)$ on $\Gamma_i$ ($i=1,2$).

\smallskip

\item Local Stokes problem depending on auxiliary interface data: find $\bu_{i,h}^\lambda(\bmu)$ and $p_{i,h}^\lambda(\bmu)$ such that
\begin{subequations}\label{eq:StokesMultiDomainLambda}
	\begin{eqnarray}
		-\nabla \cdot \bsigma(\bu_{i,h}^\lambda(\bmu), p_{i,h}^\lambda(\bmu); \bmu) = \mathbf{0} &&\quad \text{in } \Omega_i,\\
		\nabla \cdot \bu_{i,h}^\lambda(\bmu) = 0 &&\quad \text{in } \Omega_i,\\
		\bu_{i,h}^\lambda(\bmu) = \mathbf{0} &&\quad \text{on } \Gamma^D_i,\\
		\bsigma(\bu_{i,h}^\lambda(\bmu), p_{i,h}^\lambda(\bmu); \bmu) \bn = \mathbf{0} &&\quad \text{on } \Gamma^N_i,\\
		\bu_{i,h}^\lambda(\bmu) = \blambda_{i,h}(\bmu) &&\quad \text{on } \Gamma_i\,.
	\end{eqnarray}
\end{subequations}
The auxiliary trace function $\blambda_{i,h}(\bmu) \in \bY_{i}^{h}$ must be determined to guarantee that
\begin{equation*}
    \bu_{i,h}(\bmu) = \bu_{i,h}^f(\bmu) + \bu_{i,h}^\lambda(\bmu)
    \quad \text{and} \quad
    p_{i,h}(\bmu) = p_{i,h}^f(\bmu) + p_{i,h}^\lambda(\bmu)
    \quad \text{in } \Omega_i\,,
\end{equation*}
where $\bu_{i,h}(\bmu)$ and $p_{i,h}(\bmu)$ are the Galerkin approximations of the velocity and pressure in \eqref{eq:StokesMultiDomain}. (Note that, due to the definition of $\bY_{i}^{h}$, $\blambda_{i,h}(\bmu) = \mathbf{0}$ on $\overline{\Gamma}_i^D \cap \overline{\Gamma}_i$ if $\overline{\Gamma}_i^D \cap \overline{\Gamma}_i \not= \emptyset$ to avoid a possible discontinuity in the boundary condition for the velocity $\bu_{i,h}^\lambda(\bmu)$.)
\end{enumerate}

Let $\blambda_{\Omega_i,h} (\bmu) \in [H^1(\Omega_i)\cap X_{i}^r]^d$ denote a suitable continuous extension of $\blambda_{i,h} (\bmu)$ such that $\blambda_{\Omega_i,h} (\bmu)_{\vert\Gamma_i} = \blambda_{i,h}(\bmu)$ and $\blambda_{\Omega_i,h} (\bmu) = \mathbf{0}$ on $\Gamma_i^D$ for $i=1,2$.
A possible way to construct the extensions $\blambda_{\Omega_i,h}(\bmu)$ is the following. For each $k=1,\ldots,d$, let $\Lambda_i^{k,j}(\bmu)$, with $j=1,\ldots,N_{\Gamma_i}$, be the nodal values of the $k$th component $\lambda_{i,h}^k(\bmu)$ of $\blambda_{i,h}(\bmu)$ at the $N_{\Gamma_i}$ degrees of freedom (dofs) on $\Gamma_i$, and let $\bLambda_{\Gamma_i}$ be the $d \times N_{\Gamma_i}$ dimensional vector that collects such coefficients. Then, we can write $\blambda_{i,h}(\bmu) = (\lambda_{i,h}^1(\bmu),\ldots,\lambda_{i,h}^d(\bmu))$ with
\begin{equation}
	\label{eq:lambda}
	\lambda_{i,h}^k(\bmu) = \lambda_{i,h}^k(\bx;\bmu) = \sum_{j = 1}^{N_{\Gamma_i}} \Lambda_i^{k, j}(\bmu)\, \eta_{i,h}^j(\bx), \quad k = 1, \dots, d,
\end{equation}
where $\eta_{i,h}^j(\bx)$ are the basis functions of $Y_i^h$ associated with the $N_{\Gamma_i}$ nodes on $\Gamma_i$. A possible strategy to select the basis functions $\eta_{i,h}^j(\bx)$ is to consider the non-null restrictions on $\Gamma_i$ of the basis functions, say, $\varphi_{i,h}^j(\bx)$, used inside $\Omega_i$ to discretize the various components of the velocity: $\eta_{i,h}^j(\bx) = \varphi_{i,h}^j(\bx)_{\vert \Gamma_i}$. We can then express $\blambda_{\Omega_i,h}(\bmu) = (\lambda_{\Omega_i,h}^1(\bmu),\ldots,\lambda_{\Omega_i,h}^d(\bmu))$, with each component defined as
\begin{equation}\label{eq:lambdaOmega}
    \lambda_{\Omega_i,h}^k(\bmu) = \lambda_{\Omega_i,h}^k(\bx;\bmu) =
    \sum_{j = 1}^{N_{\Gamma_i}} \Lambda_i^{k, j}(\bmu) \, \varphi_{i,h}^j(\bx), \quad k = 1, \dots, d.
\end{equation}
Whilst this is the choice made in the rest of this work due to its simple implementation and non-intrusiveness with respect to the underlying finite element code, other suitable bases can be used (see \cite{Discacciati:2024:CMAME} for discussion).

Then, we can split
\begin{equation*}
    \bu_{i,h}^f(\bmu) = \bu_{i,h}^{0,f}(\bmu) + \bg_{\Omega_i}^D(\bmu)
    \quad \text{and} \quad
    \bu_{i,h}^\lambda(\bmu) = \bu_{i,h}^{0,\lambda}(\bmu) + \blambda_{\Omega_i,h}(\bmu),
\end{equation*}
with $\bu_{i,h}^{0,f}(\bmu), \bu_{i,h}^{0,\lambda}(\bmu) \in \bV_i^h$. Denoting 
\begin{equation*}
\bu_{i,h}^0(\bmu) = \bu_{i,h}^{0,f}(\bmu) +  \bu_{i,h}^{0,\lambda}(\bmu) \in \bV_i^h,    
\end{equation*}
the Galerkin approximation of the velocity $\bu_{i}(\bmu)$ in \eqref{eq:StokesMultiDomain} can finally be expressed as 
\begin{equation}\label{eq:splitting}
    \bu_{i,h}(\bmu) = \bu_{i,h}^0(\bmu) + \blambda_{\Omega_i,h}(\bmu) + \bg_{\Omega_i}^D(\bmu) \, .
\end{equation}
Therefore, using \eqref{eq:splitting}, the discrete counterpart of the interface condition \eqref{eq:StokesMultiDomain_5} can be written as
\begin{subequations}\label{eq:StokesMultiDomainWeakConditions}
\begin{eqnarray}
 \blambda_{1,h}(\bmu) - \blambda_{\Omega_2,h}(\bmu)_{\vert\Gamma_1} - \bu^0_{2,h}(\bmu)_{\vert\Gamma_1} &=& (- \bg^D_{\Omega_1}(\bmu) + \bg^D_{\Omega_2}(\bmu))_{\vert\Gamma_1} \quad \text{on } \Gamma_1, \\
 \blambda_{2,h}(\bmu) - \blambda_{\Omega_1,h}(\bmu)_{\vert\Gamma_2} - \bu^0_{1,h}(\bmu)_{\vert\Gamma_2}   &=& (\phantom{-}\bg^D_{\Omega_1}(\bmu) - \bg^D_{\Omega_2}(\bmu))_{\vert\Gamma_2} \quad \text{on } \Gamma_2.
\end{eqnarray}
\end{subequations}

For all $\bu,\bv \in [H^1(\Omega_i)]^d$ and $p,q \in L^2(\Omega_i)$, consider the bilinear forms
\begin{eqnarray*}
\mathcal{A}_i(\bu, \bv; \bmu) &=&
   \int_{\Omega_i} 2\nu(\bmu) \nabla^s \bu : \nabla^s \bv \, d\bx +
   \mathcal{A}_i^{stab}(\bu, \bv; \bmu) \,, \\
\mathcal{B}_i(q, \bv; \bmu) &=&
   - \int_{\Omega_i} q \, \nabla \cdot \bv \, d\bx +    \mathcal{B}_i^{stab}(q, \bv; \bmu)\, ,\\
\mathcal{C}_i(p, q) &=&
   \mathcal{C}_i^{stab}(p, q)\,,
\end{eqnarray*}
and the functionals
\begin{equation*}
\mathcal{F}_i(\bv; \bmu) =
   \int_{\Omega_i} \beff_i(\bmu) \cdot \bv  \, d\bx + \int_{\Gamma_i^N} \bg_i^N(\bmu) \cdot \bv  \, d\bx +
   \mathcal{F}_i^{stab}(\bv; \bmu),
\quad
\mathcal{G}_i(q;\bmu) =
   \mathcal{G}_i^{stab}(q;\bmu)\,,
\end{equation*}
where all the contributions denoted by ${}^{stab}$ are either identically null if the spaces $\bV_i^h$ and $W_i^h$ are inf-sup stable, or they introduce suitable stabilizing terms if the inf-sup condition is not satisfied, following the classical finite element theory for the Stokes problem (see, e.g., \cite{Quarteroni:1994,Donea:2003}). For example, in the case of the Galerkin Least Squares stabilization, for all $\bu,\bv \in [H^1(\Omega_i)]^d$ and $q \in L^2(\Omega_i)$, one would have
\begin{eqnarray*}
\mathcal{A}_i^{stab} (\bu,\bv;\bmu) &=& - \, \delta \sum_{K_i} h_{K_i}^2 \, \int_{K_i} (\nabla \cdot (\nu(\bmu) \nabla^s \bu)) \cdot (\nabla \cdot (\nu(\bmu) \nabla^s \bv)) \, d\bx \, , \\
\mathcal{B}_i^{stab} (q,\bv; \bmu) &=& \phantom{-} \delta \sum_{K_i} h_{K_i}^2 \, \int_{K_i} \nabla q \cdot (\nabla \cdot (\nu(\bmu) \nabla^s \bv)) \, d\bx \, , \\
\mathcal{C}_i^{stab} (p,q) &=& - \, \delta \sum_{K_i} h_{K_i}^2 \, \int_{K_i} \nabla p \cdot \nabla q \, d\bx \, , \\
\mathcal{F}_i^{stab} (\bv;\bmu) &=& \phantom{-} \delta \sum_{K_i} h_{K_i}^2\, \int_{K_i} \beff_i (\bmu) \cdot (\nabla \cdot (\nu(\bmu) \nabla^s \bv)) \, d\bx \,, \\
\mathcal{G}_i^{stab} (q;\bmu) &=& - \, \delta \sum_{K_i} h_{K_i}^2 \int_{K_i} \beff_i(\bmu) \cdot \nabla q  \, d\bx\, ,
\end{eqnarray*}
where $K_i \in \mathcal{T}_i$ is a generic element of size $h_{K_i}$ of the computational mesh, and $\delta > 0$ is a positive stabilization coefficient.

Then, the Galerkin finite element approximation of the two-domain formulation \eqref{eq:StokesMultiDomain} with interface conditions \eqref{eq:StokesMultiDomain_5} expressed as in \eqref{eq:StokesMultiDomainWeakConditions} becomes: for all $\bmu \in \mathcal{P}$, for $i=1,2$, find $(\bu_{i,h}^{0}(\bmu), p_{i,h}(\bmu)) \in \bV_i^h \times W_i^h$ and $\blambda_{i,h}(\bmu) \in \bY_i^h$ such that, for all $(\bv_{i,h}, q_{i,h}) \in \bV_i^h \times W_i^h$,
\begin{subequations}\label{eq:StokesMultiDomainGalerkin}
	\begin{eqnarray}
    \hspace*{-7mm}
		\mathcal{A}_i(\bu_{i,h}^{0}(\bmu) + \blambda_{\Omega_i,h}(\bmu), \bv_{i,h}; \bmu) + \mathcal{B}_i(p_{i,h}(\bmu), \bv_{i,h}; \bmu) &=& \mathcal{F}_i(\bv_{i,h};\bmu) - \mathcal{A}_i (\bg_{\Omega_i}^D(\bmu), \bv_{i,h} ; \bmu),\label{eq:StokesMultiDomainGalerkin_1}\\[5pt]
        \mathcal{B}_i(q_{i,h},\bu_{i,h}^0(\bmu) + \blambda_{\Omega_i,h}(\bmu); \bmu) + \mathcal{C}_i (p_{i,h}(\bmu),q_{i,h}) &=& \mathcal{G}_i (q_{i,h};\bmu) - \mathcal{B}_i (q_{i,h},\bg_{\Omega_i}^D(\bmu); \bmu)\,, \label{eq:StokesMultiDomainGalerkin_2}\\[5pt]
        \blambda_{1,h}(\bmu) - \blambda_{\Omega_2,h}(\bmu)_{\vert\Gamma_1} - \bu^0_{2,h}(\bmu)_{\vert\Gamma_1} &=& (- \bg^D_{\Omega_1}(\bmu) + \bg^D_{\Omega_2}(\bmu))_{\vert\Gamma_1} \quad \text{on } \Gamma_1, \label{eq:StokesMultiDomainGalerkin_3}\\[5pt]
        \blambda_{2,h}(\bmu) - \blambda_{\Omega_1,h}(\bmu)_{\vert\Gamma_2} - \bu^0_{1,h}(\bmu)_{\vert\Gamma_2}   &=& (\phantom{-}\bg^D_{\Omega_1}(\bmu) - \bg^D_{\Omega_2}(\bmu))_{\vert\Gamma_2} \quad \text{on } \Gamma_2. \label{eq:StokesMultiDomainGalerkin_4}
	\end{eqnarray}
\end{subequations}

\medskip

To write the algebraic form of \eqref{eq:StokesMultiDomainGalerkin}, we use the following notation. If $\mat{M}_{\all\all}^i$ denotes a generic finite element matrix associated with the dofs in $\Omega_i$, we replace the sub-index $\all$ by $I$ (respectively, $\Gamma_i$) to denote that only the rows/columns associated with the dofs internal to $\Omega_i$ excluding $\Gamma_i$ (respectively, the dofs on $\Gamma_i$) are retained, with the index for rows predecing the one for columns. Moreover, let $\mat{I}_{\Gamma_i\Gamma_i}$ be the identity matrix associated with the dofs on $\Gamma_i$, and $\mat{R}_{\Omega_i\to\Gamma_j}$ ($i,j=1,2$, $i \not= j$) be the restriction matrix that, given any vector of nodal values in $\Omega_i$, returns the vector of nodal values on the interface $\Gamma_j$ inside $\Omega_i$. 
Finally, we denote by $\mat{E}_{\Gamma_i\to\Omega_i} \bLambda_{\Gamma_i} \in [X_i^r]^d$, with $\mat{E}_{\Gamma_i\to\Omega_i} \bLambda_{\Gamma_i} = \mathbf{0}$ on $\Gamma_i^D$, a suitable algebraic extension operator of the interface nodal values $\bLambda_{\Gamma_i}$ that corresponds, e.g., to \eqref{eq:lambdaOmega} at the algebraic level.

Then, with additional self-explanatory notation, the linear system corresponding to the Galerkin approximation \eqref{eq:StokesMultiDomainGalerkin} becomes
\begin{equation}\label{eq:StokesMultiDomainAlgebraic}
 \begin{pmatrix}
 \mat{A}^1_{II} & (\mat{B}^1_{\all I})^T & \mat{0} & \mat{0} & \mat{A}^1_{I\Gamma_1} & \mat{0} \\[2pt]
 \mat{B}^1_{\all I} & \mat{C}^1_{\all \all} & \mat{0} & \mat{0} & \mat{B}^1_{\all \Gamma_1} & \mat{0} \\[2pt]
 \mat{0} & \mat{0} & \mat{A}^2_{II} & (\mat{B}^2_{\all I})^T & \mat{0} & \mat{A}^2_{I\Gamma_2} \\[2pt]
 \mat{0} & \mat{0} & \mat{B}^2_{\all I} & \mat{C}^2_{\all \all} & \mat{0} & \mat{B}^2_{\all \Gamma_2} \\[2pt]
 \mat{0} & \mat{0} & -\mat{R}_{\Omega_2 \to \Gamma_1} & \mat{0} & \mat{I}_{\Gamma_1\Gamma_1} & -\mat{R}_{\Omega_2 \to \Gamma_1} \mat{E}_{\Gamma_2 \to \Omega_2} \\[2pt]
 -\mat{R}_{\Omega_1\to\Gamma_2} & \mat{0} & \mat{0} & \mat{0} & -\mat{R}_{\Omega_1 \to \Gamma_2} \mat{E}_{\Gamma_1 \to \Omega_1} & \mat{I}_{\Gamma_2\Gamma_2}
 \end{pmatrix}
 \begin{pmatrix}
     \mathbf{u}^1_I \\[2pt]
     \mathbf{p}^1_\all \\[2pt]
     \mathbf{u}^2_I \\[2pt]
     \mathbf{p}^2_\all \\[2pt]
     \boldsymbol{\Lambda}_{\Gamma_1} \\[2pt]
     \boldsymbol{\Lambda}_{\Gamma_2}
 \end{pmatrix}
 =
 \begin{pmatrix}
     \mathbf{f}^1_I \\[2pt]
     \mathbf{g}^1_\all \\[2pt]
     \mathbf{f}^2_I \\[2pt]
     \mathbf{g}^2_\all \\[2pt]
     \mathbf{g}_{\Gamma_1} \\[2pt]
     \mathbf{g}_{\Gamma_2}
 \end{pmatrix}.
\end{equation}

Finally, by static condensation, we can equivalently rewrite the algebraic problem \eqref{eq:StokesMultiDomainAlgebraic} in terms of the interface variables only. This yields the interface system
\begin{equation}\label{eq:StokesMultiDomainAlgebraicInterface}
\begin{array}{l}
\left(
 \begin{pmatrix}
  \mat{I}_{\Gamma_1\Gamma_1} & -\mat{R}_{\Omega_2 \to \Gamma_1} \mat{E}_{\Gamma_2 \to \Omega_2} \\
  -\mat{R}_{\Omega_1 \to \Gamma_2} \mat{E}_{\Gamma_1 \to \Omega_1} & \mat{I}_{\Gamma_2\Gamma_2}
 \end{pmatrix}
\right. \\
\qquad -
 \begin{pmatrix}
  \mat{0} & \mat{0} \\
  \mat{R}_{\Omega_1\to\Gamma_2} & \mat{0}
 \end{pmatrix}
 \begin{pmatrix}
  \mat{A}^1_{II} & (\mat{B}^1_{\all I})^T \\
  \mat{B}^1_{\all I} & \mat{C}^1_{\all \all} 
 \end{pmatrix}^{-1}
 \begin{pmatrix}
  -\mat{A}^1_{I\Gamma_1} & \mat{0} \\
  -\mat{B}^1_{\all \Gamma_1} & \mat{0} 
 \end{pmatrix} \\
\qquad -
\left.
 \begin{pmatrix}
  \mat{R}_{\Omega_2 \to \Gamma_1} & \mat{0} \\
  \mat{0} & \mat{0} 
 \end{pmatrix}
 \begin{pmatrix}
  \mat{A}^2_{II} & (\mat{B}^2_{\all I})^T \\
  \mat{B}^2_{\all I} & \mat{C}^2_{\all \all}
 \end{pmatrix}^{-1}
 \begin{pmatrix}
  \mat{0} & -\mat{A}^2_{I\Gamma_2} \\
  \mat{0} & -\mat{B}^2_{\all \Gamma_2}
 \end{pmatrix}
\right)
\begin{pmatrix}
  \boldsymbol{\Lambda}_{\Gamma_1} \\
  \boldsymbol{\Lambda}_{\Gamma_2}
\end{pmatrix} \\
= 
 \begin{pmatrix}
  \mathbf{g}_{\Gamma_1} \\
  \mathbf{g}_{\Gamma_2}
 \end{pmatrix}
+
 \begin{pmatrix}
  \mat{0} & \mat{0} \\
  \mat{R}_{\Omega_1\to\Gamma_2} & \mat{0}
 \end{pmatrix}
 \begin{pmatrix}
  \mat{A}^1_{II} & (\mat{B}^1_{\all I})^T \\
  \mat{B}^1_{\all I} & \mat{C}^1_{\all \all} 
 \end{pmatrix}^{-1}
 \begin{pmatrix}
  \mathbf{f}^1_I \\
  \mathbf{g}^1_\all \\
 \end{pmatrix} \\
\phantom{= 
 \begin{pmatrix}
  \mathbf{g}_{\Gamma_1} \\
  \mathbf{g}_{\Gamma_2}
 \end{pmatrix}}
+
 \begin{pmatrix}
  \mat{R}_{\Omega_2 \to \Gamma_1} & \mat{0} \\
  \mat{0} & \mat{0} 
 \end{pmatrix}
 \begin{pmatrix}
  \mat{A}^2_{II} & (\mat{B}^2_{\all I})^T \\
  \mat{B}^2_{\all I} & \mat{C}^2_{\all \all}
 \end{pmatrix}^{-1}
 \begin{pmatrix}
  \mathbf{f}^2_I \\
  \mathbf{g}^2_\all
\end{pmatrix}.
\end{array}
\end{equation}
Observe that, for $i=1,2$, applying the inverse matrix
\begin{equation*}
 \begin{pmatrix}
  \mat{A}^i_{II} & (\mat{B}^i_{\all I})^T \\
  \mat{B}^i_{\all I} & \mat{C}^i_{\all \all}
 \end{pmatrix}^{-1}
\end{equation*}
implies solving a Stokes problem in the subdomain $\Omega_i$. Therefore, to construct the right-hand side of system \eqref{eq:StokesMultiDomainAlgebraicInterface}, one must solve the two local Stokes problems
\begin{equation}\label{eq:localStokesAlgebraicData}
 \begin{pmatrix}
  \mat{A}^i_{II} & (\mat{B}^i_{\all I})^T \\
  \mat{B}^i_{\all I} & \mat{C}^i_{\all \all}
 \end{pmatrix}^{-1}
 \begin{pmatrix}
  \mathbf{f}^i_I \\
  \mathbf{g}^i_\all
\end{pmatrix},
\quad i=1,2,
\end{equation}
that depend on the problem data but that are independent of the interface variables $\bLambda_{\Gamma_i}$. Notice that these problems are the algebraic counterpart of \eqref{eq:StokesMultiDomainF}.

System \eqref{eq:StokesMultiDomainAlgebraicInterface} can be solved by an iterative matrix-free Krylov method (e.g., GMRES) and, at each iteration $n$, for $i=1,2$, the local problems
\begin{equation}\label{eq:localStokesAlgebraicLambda}
 \begin{pmatrix}
  \mat{A}^i_{II} & (\mat{B}^i_{\all I})^T \\
  \mat{B}^i_{\all I} & \mat{C}^i_{\all \all}
 \end{pmatrix}^{-1}
 \begin{pmatrix}
  -\mat{A}^i_{I\Gamma_i} \\
  -\mat{B}^i_{\all \Gamma_i}
 \end{pmatrix}
 \boldsymbol{\Lambda}_{\Gamma_i}^{(n)}
\end{equation}
must be solved. These are the algebraic counterpart of the local Stokes problems \eqref{eq:StokesMultiDomainLambda} with homogeneous data and homogeneous boundary conditions, except for the Dirichlet boundary conditions on $\Gamma_i$ that impose that the discretized velocity is equal to the nodal values $\boldsymbol{\Lambda}_{\Gamma_i}^{(n)}$.

Solving problems \eqref{eq:localStokesAlgebraicData} and \eqref{eq:localStokesAlgebraicLambda} with a full-order solver (e.g., the finite element method) makes each iteration of the Krylov method computationally expensive. In the parametric setting that we are considering here, the computational cost becomes even higher because, whenever a new set of parameters, say, $\bar{\bmu} \in \mathcal{P}$, is considered, the whole procedure must be executed from the beginning, possibly starting from the assembly of the matrices since these might depend on the parameters themselves.

Therefore, following the approach proposed in \cite{Discacciati:2024:CMAME,Discacciati:2024:DD28}, we reformulate the algorithm to solve the interface problem \eqref{eq:StokesMultiDomainAlgebraicInterface} by identifying an offline phase and an online phase to efficiently handle the presence of the physical parameters $\bmu \in \mathcal{P}$. More precisely, in the offline phase, the local parametric problems~\eqref{eq:localStokesAlgebraicData} and~\eqref{eq:localStokesAlgebraicLambda} are solved by PGD to obtain local surrogate models that, in the case of \eqref{eq:localStokesAlgebraicLambda}, must incorporate arbitrary traces $\blambda_i$ on the interfaces. The procedure needed to construct such local surrogate models is detailed in Sect.~\ref{sect:StokesSurrogateModel}. In the online phase, the interface system \eqref{eq:StokesMultiDomainAlgebraicInterface} is solved by GMRES for any fixed set of parametric values $\bar{\bmu} \in \mathcal{P}$. However, at this stage, the solution of the local problems~\eqref{eq:localStokesAlgebraicData} and~\eqref{eq:localStokesAlgebraicLambda} is replaced by the \textit{evaluation} of the precomputed local surrogate models, thus allowing a real-time execution of the coupling algorithm. Details are provided in Sect.~\ref{sect:StokesAScoupling}.

\smallskip

We conclude this section by remarking that other possible multi-domain formulations than \eqref{eq:StokesMultiDomain} can be considered for the Stokes problem \eqref{eq:StokesGlobal} as discussed in \cite{Discacciati:2013:Stokes}. Indeed, one could impose the continuity of the normal component of the Cauchy stress tensor on $\Gamma_1\cup\Gamma_2$ instead of \eqref{eq:StokesMultiDomain_5}, or the continuity of the velocity, say, on $\Gamma_1$, and of the normal stress on $\Gamma_2$. The associated interface system would then be modified accordingly. As pointed out in \cite{Discacciati:2013:Stokes}, imposing the continuity of the normal stress on at least one of the interfaces generally results in a lower number of iterations of the coupling algorithm, because this approach leads to a better control not only of the velocity but also of the pressure on the interfaces and across the whole overlapping region. However, the choice of considering here only the continuity of the velocity traces on the interfaces is due to the simplicity of handling parametric boundary conditions of Dirichlet type in the context of a continuous Galerkin finite element approximation of the Stokes problem. Indeed, as we will see in Sect.~\ref{sect:StokesSurrogateModel}, we can easily parametrise the interface conditions \eqref{eq:StokesMultiDomainGalerkin_3} and \eqref{eq:StokesMultiDomainGalerkin_4} in the offline phase, and efficiently identify and exchange information about the velocity on the interfaces in the online phase without having to solve any interface problem. This would not be the case if Neumann boundary conditions -- which would be imposed weakly on the interfaces -- were considered. A more convenient way to handle natural boundary conditions in the computational framework developed in this work would be to use different discretization techniques such as, e.g., discontinuous Galerkin methods. This will be the object of future work.

\subsection{Local surrogate models using proper generalized decomposition for the Stokes problem}
\label{sect:StokesSurrogateModel}

In this section, we focus on the construction of PGD-based local surrogate models for the Stokes problems \eqref{eq:localStokesAlgebraicData} and \eqref{eq:localStokesAlgebraicLambda} in the offline phase.

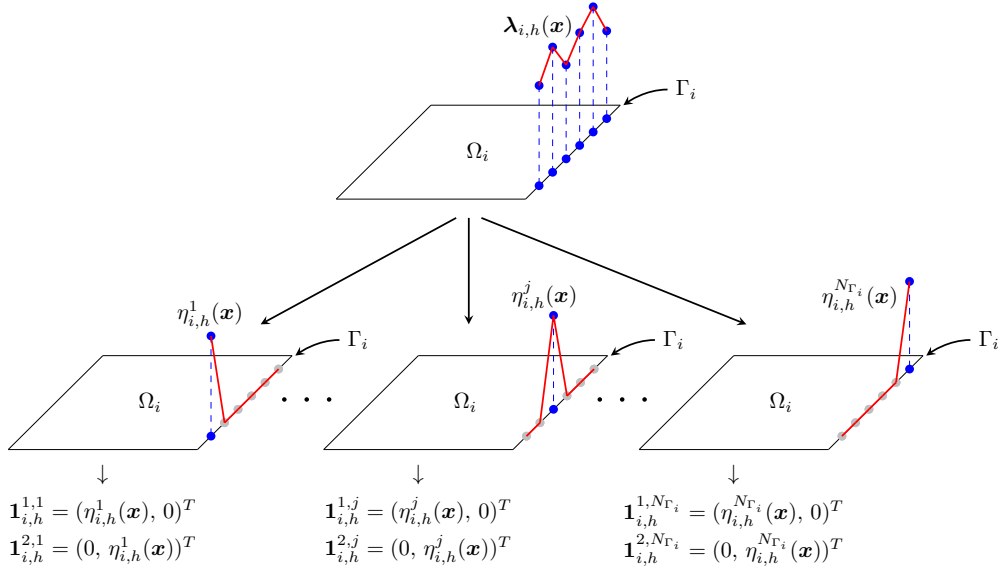
\begin{figure}[hbt]
\begin{center}
\resizebox{0.9\textwidth}{!}{
	\begin{tikzpicture}
		\newcommand{\intNodes}[3]{
				\draw[dashed, #3] (\longSide + #1 * \delta / 7, \ang * #1 * \delta / 7 + #2) -- (\longSide + #1 * \delta / 7, #1 * \ang * \delta / 7 );
				\fill[#3] (\longSide + #1 * \delta / 7, \ang * #1 * \delta / 7 + #2) circle (0.07 cm);
				\fill[#3] (\longSide + #1 * \delta / 7, #1 * \ang * \delta / 7 ) circle (0.07 cm);
			}
			
		\newcommand{\intFun}[3]{
				\draw[thick, red] (\longSide + #1 * \delta / 7, #1 * \ang * \delta / 7 + #2 ) -- (\longSide + #1 * \delta / 7 + \delta / 7, #1 * \ang * \delta / 7 +  \ang * \delta / 7 + #3);
			}
			
			\def\ang{1}
			\def\delta{1.5}
			\def\longSide{3}
			
			\begin{scope}[shift={(0.2, 0)}]
				\draw (0, 0) -- (\delta, \delta*\ang);
				\draw (0, 0) -- (\longSide, 0);
				\draw(\delta, \delta*\ang) -- (\delta + \longSide, \delta*\ang);
				\draw(\longSide, 0) -- (\delta + \longSide, \delta*\ang);
				
				\node at (\longSide / 2 + \delta / 2, \delta * \ang / 2) {$\Omega_i$};

                \node at (\longSide / 2 + 1.7, \delta * \ang / 2 + 2) {$\blambda_{i,h}(\bx)$};
				
				\draw[-stealth, thick] (\longSide + \delta + 0.75, \delta * \ang + 0.25) node[anchor=west] {$\Gamma_i$} parabola (\longSide + \delta + 0.05, \delta * \ang);
				
				\intNodes{1}{1.6}{blue};
				\intNodes{2}{2}{blue};
				\intNodes{3}{1.5}{blue};	
				\intNodes{4}{1.8}{blue};
				\intNodes{5}{2}{blue};
				\intNodes{6}{1.4}{blue};
				
				\intFun{1}{1.6}{2};
				\intFun{2}{2}{1.5};
				\intFun{3}{1.5}{1.8};
				\intFun{4}{1.8}{2};
				\intFun{5}{2}{1.4};
			\end{scope}
			
			\begin{scope}[shift={(-5, -4)}]
				\draw (0, 0) -- (\delta, \delta*\ang);
				\draw (0, 0) -- (\longSide, 0);
				\draw(\delta, \delta*\ang) -- (\delta + \longSide, \delta*\ang);
				\draw(\longSide, 0) -- (\delta + \longSide, \delta*\ang);
				
				\node at (\longSide / 2 + \delta / 2, \delta * \ang / 2) {$\Omega_i$};

                \node at (\longSide / 2 + 1.7, \delta * \ang / 2 + 1.4) {$\eta_{i,h}^1(\bx)$};
				
				\draw[-stealth, thick] (\longSide + \delta + 0.75, \delta * \ang + 0.25) node[anchor=west] {$\Gamma_i$} parabola (\longSide + \delta + 0.05, \delta * \ang);
				
				\intNodes{1}{1.6}{blue};
				\intNodes{2}{0}{lightgray};
				\intNodes{3}{0}{lightgray};
				\intNodes{4}{0}{lightgray};
				\intNodes{5}{0}{lightgray};
				\intNodes{6}{0}{lightgray};

				\intFun{1}{1.6}{0};
				\intFun{2}{0}{0};
				\intFun{3}{0}{0};
				\intFun{4}{0}{0};
				\intFun{5}{0}{0};
			\end{scope}
		
			\begin{scope}[shift={(0, -4)}]
				\draw (0, 0) -- (\delta, \delta*\ang);
				\draw (0, 0) -- (\longSide, 0);
				\draw(\delta, \delta*\ang) -- (\delta + \longSide, \delta*\ang);
				\draw(\longSide, 0) -- (\delta + \longSide, \delta*\ang);
				
				\node at (\longSide / 2 + \delta / 2, \delta * \ang / 2) {$\Omega_i$};

                \node at (\longSide / 2 + 2, \delta * \ang / 2 + 1.7) {$\eta_{i,h}^j(\bx)$};
				
				\draw[-stealth, thick] (\longSide + \delta + 0.75, \delta * \ang + 0.25) node[anchor=west] {$\Gamma_i$} parabola (\longSide + \delta + 0.05, \delta * \ang);
				
				\intNodes{1}{0}{lightgray};
				\intNodes{2}{0}{lightgray};
				
				\intNodes{3}{1.5}{blue};
				
				\intNodes{4}{0}{lightgray};
				\intNodes{5}{0}{lightgray};
				\intNodes{6}{0}{lightgray};

				\intFun{1}{0}{0};
				\intFun{2}{0}{1.5};
				\intFun{3}{1.5}{0};
				\intFun{4}{0}{0};
				\intFun{5}{0}{0};
			\end{scope}
		
			\begin{scope}[shift={(5, -4)}]
				\draw (0, 0) -- (\delta, \delta*\ang);
				\draw (0, 0) -- (\longSide, 0);
				\draw(\delta, \delta*\ang) -- (\delta + \longSide, \delta*\ang);
				\draw(\longSide, 0) -- (\delta + \longSide, \delta*\ang);
				
				\node at (\longSide / 2 + \delta / 2, \delta * \ang / 2) {$\Omega_i$};

                \node at (\longSide / 2 + 2, \delta * \ang / 2 + 1.7) {$\eta_{i,h}^{N_{\Gamma_i}}(\bx)$};
				
				\draw[-stealth, thick] (\longSide + \delta + 0.75, \delta * \ang + 0.25) node[anchor=west] {$\Gamma_i$} parabola (\longSide + \delta + 0.05, \delta * \ang);
				
				\intNodes{1}{0}{lightgray};
				\intNodes{2}{0}{lightgray};	
				\intNodes{3}{0}{lightgray};
				\intNodes{4}{0}{lightgray};
				\intNodes{5}{0}{lightgray};
				\intNodes{6}{1.4}{blue};

				\intFun{1}{0}{0};
				\intFun{2}{0}{0};
				\intFun{3}{0}{0};
				\intFun{4}{0}{0};
				\intFun{5}{0}{1.4};
			\end{scope}
			
			\draw[- stealth, thick] (2.1, -0.3) -- (-1 , -2);
			\draw[- stealth, thick] (2.3, -0.3) -- (2.3 , -2);
			\draw[- stealth, thick] (2.5, -0.3) -- (6.7 , -2);
			
			\node at (-0.2, -3.2) {\Huge \ldots};
			\node at (4.8, -3.2) {\Huge \ldots};
			
			\node at (-3.5, -4.4) {$\downarrow$};
			\node at ( 1.5, -4.4) {$\downarrow$};
			\node at ( 6.5, -4.4) {$\downarrow$};
			
			\node at (-3.5, -5.0) {$\bm{1}_{i,h}^{1,1} = (\eta_{i,h}^1(\bx), \, 0 )^T$};
            \node at (-3.5, -5.6) {$\bm{1}_{i,h}^{2,1} = (0, \, \eta_{i,h}^1(\bx))^T$};
            \node at ( 1.5, -5.0) {$\bm{1}_{i,h}^{1,j} = (\eta_{i,h}^j(\bx), \, 0 )^T$};
            \node at ( 1.5, -5.6) {$\bm{1}_{i,h}^{2,j} = (0, \, \eta_{i,h}^j(\bx))^T$};
            
            \node at ( 6.5, -5.0) {$\bm{1}_{i,h}^{1,N_{\Gamma_i}} = (\eta_{i,h}^{N_{\Gamma_i}}(\bx), \, 0 )^T$};
            \node at ( 6.5, -5.6) {$\bm{1}_{i,h}^{2,N_{\Gamma_i}} = (0, \, \eta_{i,h}^{N_{\Gamma_i}}(\bx))^T$};
			
		\end{tikzpicture}
        }
\end{center}
		\caption{Sketch illustrating the basis functions $\eta_{i,h}^j(\bx)$ and the corresponding discrete vectorial trace function $\bm{1}_{i,h}^{k, j} \in \bY_i^h$ in a 2D setting.}
		\label{fig:activeBdryNodes}
\end{figure}

We begin by introducing the discrete vectorial trace function $\bm{1}_{i,h}^{k, j} \in \bY_i^h$ whose $d$ components are identically null on $\Gamma_i$, except for its $k$th component that is equal to the basis function $\eta_{i,h}^j(\bx)$ (see \eqref{eq:lambda}) on $\Gamma_i$. (See Fig.~\ref{fig:activeBdryNodes}.) 

Moreover, let $\bm{1}_{\Omega_i,h}^{k,j} \in [X_{i}^r]^d$ be a suitable continuous extension of $\bm{1}_{i,h}^{k,j}$ in $\Omega_i$. Then, consider the local subproblems: for all $\bmu \in \mathcal{P}$, for $i=1,2$, $k=1,\ldots,d$ and $j=1,\ldots,N_{\Gamma_i}$, find $({\bu}_{i,h}^{k,j}(\bmu), {p}_{i,h}^{k,j}(\bmu)) \in \bV_i^h \times W_i^h$ such that, for all $(\bv_{i,h}, q_{i,h}) \in \bV_i^h \times W_i^h$,
\begin{subequations}\label{eq:StokesMultiDomainGalerkinUnity}
	\begin{eqnarray}
		\mathcal{A}_i({\bu}_{i,h}^{k,j}(\bmu), \bv_{i,h}; \bmu) + \mathcal{B}_i({p}_{i,h}^{k,j}(\bmu), \bv_{i,h}; \bmu) &=& - \mathcal{A}_i(\bm{1}_{\Omega_i,h}^{k,j}, \bv_{i,h}; \bmu)\,,\\
        \mathcal{B}_i(q_{i,h},{\bu}_{i,h}^{k,j}(\bmu); \bmu) + \mathcal{C}_i ({p}_{i,h}^{k,j}(\bmu),q_{i,h}) &=& - \mathcal{B}_i(q_{i,h},\bm{1}_{\Omega_i,h}^{k,j}; \bmu)\,.
	\end{eqnarray}
\end{subequations}
Using the expansion \eqref{eq:lambda} and exploiting the linearity of the Stokes problem, by construction there holds
\begin{equation}\label{eq:upPGDLambda}
    \bu^{0,\lambda}_{i,h}(\bmu) = \sum_{k = 1}^d\sum_{j = 1}^{N_{\Gamma_i}} \Lambda_i^{k, j}(\bmu)\, {\bu}_{i,h}^{k,j}(\bmu),
    \qquad
	p^{\lambda}_{i,h}(\bmu) = \sum_{k = 1}^d\sum_{j = 1}^{N_{\Gamma_i}} \Lambda_i^{k, j}(\bmu)\, {p}_{i,h}^{k,j}(\bmu)\,.
\end{equation}
Notice that if $\bLambda_{\Gamma_i}$ ($i = 1, 2$) satisfies the interface system \eqref{eq:StokesMultiDomainAlgebraicInterface}, the solution of \eqref{eq:StokesMultiDomainGalerkin} can be written as 
\begin{equation*}
    \bu_{i,h}^0(\bmu) = \bu^{0,f}_{i,h}(\bmu) + \bu^{0,\lambda}_{i,h}(\bmu),
    \qquad
    p_{i,h}(\bmu) = p^f_{i,h}(\bmu) + p^\lambda_{i,h}(\bmu).
\end{equation*}

Now, we can proceed to formulate the problems that must be solved to compute the local surrogate models for Stokes using PGD. As usual when constructing PGD approximations, we assume that problem data have a separable representation -- be that their given form, or the one that can be constructed numerically, see, e.g., \cite{DM-MZH:15}. Therefore, we can write
\begin{equation*}
	\begin{array}{c}
		\nu(\bmu) = \displaystyle\sum_{\ell=1}^{N_{\nu}} \xi_{\nu}^\ell(\bmu)b_{\nu}^\ell(\bx) \,, \qquad
		\beff_i(\bmu) = \displaystyle\sum_{\ell=1}^{N_f} \xi_{i, f}^\ell(\bmu) \bb_{i, f}^\ell(\bx) \,, \\
		\bg_i^N(\bmu) = \displaystyle\sum_{\ell=1}^{N_N}  \xi_{i, N}^\ell(\bmu) \bb_{i, N}^\ell(\bx) \,, \qquad
		\bg_{\Omega_i}^D(\bmu) = \displaystyle\sum_{\ell=1}^{N_D}  \xi_{i,D}^\ell(\bmu) \bb_{i,D}^\ell(\bx) \,,
	\end{array}
\end{equation*}
with the parametric modes assumed to be the product of scalar functions depending on a single parameter, e.g.,
\begin{equation*}
	\xi_{\nu}^\ell(\bmu) = \prod_{k = 1}^{N_p} \xi_{\nu}^{\ell, k}(\mu_k).
\end{equation*}

Consider the local problem \eqref{eq:StokesMultiDomainF} whose Galerkin approximation reads: for all $\bmu \in \mathcal{P}$ and for $i=1,2$, find $(\bu_{i,h}^{0,f}(\bmu), p_{i,h}^f(\bmu)) \in \bV_i^h \times W_i^h$ such that, for all $(\bv_{i,h}, q_{i,h}) \in \bV_i^h \times W_i^h$,
\begin{subequations}\label{eq:StokesLocalGalerkinData}
	\begin{eqnarray}
		\mathcal{A}_i(\bu_{i,h}^{0,f}(\bmu), \bv_{i,h}; \bmu) + \mathcal{B}_i(p_{i,h}^f(\bmu), \bv_{i,h}; \bmu) &=& \mathcal{F}_i(\bv_{i,h};\bmu) - \mathcal{A}_i (\bg_{\Omega_i}^D(\bmu), \bv_{i,h} ; \bmu)\,,\\
        \mathcal{B}_i(q_{i,h},\bu_{i,h}^{0,f}(\bmu); \bmu) + \mathcal{C}_i (p_{i,h}^f(\bmu),q_{i,h}) &=& \mathcal{G}_i (q_{i,h};\bmu) - \mathcal{B}_i (q_{i,h},\bg_{\Omega_i}^D(\bmu);\bmu)\,.
	\end{eqnarray}
\end{subequations}
In the PGD context, the contribution of the Dirichlet boundary condition due to $\bg_{\Omega_i}^D(\bmu)$ is accounted for by introducing ad-hoc smooth modes. Then, following \cite{Diez:2017:CMAME}, the velocity and the pressure are approximated as
\begin{equation}\label{eq:SepSolnVelPress}
    \bu^{0,f}_{i,h}(\bmu) \approx \bu^{0,f}_{i,\pgd}(\bmu) = \displaystyle\sum_{m = 1}^{M_u} {\bU}_{i}^m(\bx)\,{\phi}_{i}^m(\bmu) \,,
    \quad
    {p}^f_{i,h}(\bmu) \approx {p}^f_{i,\pgd}(\bmu) = \displaystyle\sum_{m = 1}^{M_p} {P}_{i}^m(\bx)\,{\phi}_{i}^m(\bmu) \,.
\end{equation}
Here, ${\bU}_{i}^m(\bx)$ and ${P}_{i}^m(\bx)$ are the $m$th spatial modes that are discretized by the Galerkin finite element method, while ${\phi}_{i}^m(\bmu)$ are the corresponding parametric modes that, following \cite{Diez:2017:CMAME}, are the same for velocity and pressure, and are discretized by pointwise collocation.

Assuming an affine parameter dependence for the bilinear forms and for the linear functionals in \eqref{eq:StokesLocalGalerkinData}, for all $\bu,\bv \in [H^1(\Omega_i)]^d$, $p,q \in L^2(\Omega_i)$ and $\bmu \in \mathcal{P}$, we define
\begin{eqnarray*}
\mathcal{A}_i^{\pgd}(\bu, \bv; \bmu) &=&
   \mathcal{A}_i^{h,\pgd}(\bu, \bv; \bmu) +
   \mathcal{A}_i^{stab,\pgd}(\bu, \bv; \bmu), \\
\mathcal{B}_i^{\pgd}(q, \bv; \bmu) &=&
   - \int_{\Omega_i} q \, \nabla \cdot \bv \, d\bx +
   \mathcal{B}_i^{stab,\pgd}(q, \bv; \bmu), \\
\mathcal{F}_i^{\pgd}(\bv; \bmu) &=&
   \mathcal{F}_i^{h,\pgd}(\bv; \bmu) +
   \mathcal{F}_i^{stab,\pgd}(\bv; \bmu), \\
\mathcal{G}_i^{\pgd}(q;\bmu) &=&
   \mathcal{G}_i^{stab,\pgd}(q;\bmu),
\end{eqnarray*}
where
\begin{eqnarray*}
 \mathcal{A}_i^{h,\pgd}(\bu, \bv; \bmu) &=&
  \sum_{\ell=1}^{N_{\nu}} \xi_{\nu}^\ell(\bmu) \int_{\Omega_i} 2 \, b_{\nu}^\ell(\bx)\, \nabla^s \bu : \nabla^s \bv \, d\bx\,,\\
 \mathcal{A}_i^{stab,\pgd} (\bu,\bv;\bmu) &=&
  - \delta \sum_{\ell=1}^{N_{\nu}} \sum_{m=1}^{N_{\nu}}\xi_{\nu}^\ell(\bmu) \, \xi_{\nu}^m(\bmu) \left( \sum_{K_i} h_{K_i}^2 \, \int_{K_i} (\nabla \cdot (b_{\nu}^\ell(\bx)\, \nabla^s \bu)) \cdot (\nabla \cdot (b_{\nu}^m(\bx) \, \nabla^s \bv)) \, d\bx \right) , \\
 \mathcal{B}_i^{stab,\pgd} (q,\bv; \bmu) &=& \phantom{-}\delta \sum_{\ell=1}^{N_\nu} \xi_\nu^\ell(\bmu) \left( \sum_{K_i} h_{K_i}^2 \, \int_{K_i} \nabla q \cdot (\nabla \cdot (b_\nu^\ell (\bx) \nabla^s \bv)) \, d\bx \right) , \\
 \mathcal{F}_i^{h,\pgd} (\bv; \bmu) &=& \displaystyle \sum_{\ell=1}^{N_f} \xi_{i, f}^\ell(\bmu) \int_{\Omega_i} \bb_{i, f}^\ell(\bx) \cdot \bv \, d\bx + \displaystyle\sum_{\ell=1}^{N_N}  \xi_{i, N}^\ell(\bmu) \int_{\Gamma_i^N} \bb_{i, N}^\ell(\bx) \cdot \bv \, d\bx \,, \\
 \mathcal{F}_i^{stab,\pgd} (\bv;\bmu) &=& \phantom{-}\delta \sum_{\ell=1}^{N_f} \sum_{m=1}^{N_\nu} \xi_{i, f}^\ell(\bmu) \xi_\nu^m (\bmu) \left( \sum_{K_i} h_{K_i}^2\, \int_{K_i} \bb_{i, f}^\ell(\bx) \cdot (\nabla \cdot (\bb_\nu^m(\bx) \nabla^s \bv)) \, d\bx \right), \\
 \mathcal{G}_i^{stab,\pgd} (q; \bmu) &=& - \delta  \sum_{\ell=1}^{N_f} \xi_{i, f}^\ell(\bmu) \left( \sum_{K_i} h_{K_i}^2 \int_{K_i} \bb_{i, f}^\ell(\bx) \cdot \nabla q \, d\bx \right) .
\end{eqnarray*}

Then, the PGD velocity ${\bu}^{0,f}_{i,\pgd}(\bmu)$ and pressure ${p}^f_{i,\pgd}(\bmu)$ are computed by solving the parametric problem
\begin{subequations}\label{eq:StokesLocalPGDData}
	\begin{eqnarray}
    \hspace*{-8mm}
		\mathcal{A}_i^{\pgd}({\bu}^{0,f}_{i,\pgd}(\bmu), \bv_{i,h}; \bmu) &+& \mathcal{B}_i^{\pgd}({p}^f_{i,\pgd}(\bmu), \bv_{i,h}; \bmu) \nonumber \\
        &=& \mathcal{F}_i^{\pgd}(\bv_{i,h};\bmu) - \mathcal{A}_i^{\pgd} (\sum_{\ell=1}^{N_D}  \xi_{i,D}^\ell(\bmu) \bb_{i,D}^\ell(\bx), \bv_{i,h} ; \bmu)\,,\\
  \hspace*{-8mm}
        \mathcal{B}_i^{\pgd}(q_{i,h},{\bu}^{0,f}_{i,\pgd}(\bmu); \bmu) &+& \mathcal{C}_i ({p}^f_{i,\pgd}(\bmu),q_{i,h}) \nonumber \\
        &=& \mathcal{G}_i^{\pgd} (q_{i,h};\bmu) - \mathcal{B}_i^{\pgd} (q_{i,h},\sum_{\ell=1}^{N_D}  \xi_{i,D}^\ell(\bmu) \bb_{i,D}^\ell(\bx); \bmu),
	\end{eqnarray}
\end{subequations}
for all $(\bv_{i,h}, q_{i,h}) \in \bV_i^h \times W_i^h$ and for all $\bmu \in \mathcal{P}$.

\smallskip

An analogous procedure permits to obtain the PGD problems associated with \eqref{eq:StokesMultiDomainGalerkinUnity}. More precisely, we compute the local surrogate models $\bu^{k,j}_{i,\pgd}(\bmu)$ and $p^{k,j}_{i,\pgd}(\bmu)$ by solving the local problem
\begin{subequations}\label{eq:StokesLocalPGDUnity}
	\begin{eqnarray}
		\mathcal{A}_i^{\pgd}(\bu^{k,j}_{i,\pgd}(\bmu), \bv_{i,h}; \bmu) + \mathcal{B}_i^{\pgd}(p^{k,j}_{i,\pgd}(\bmu), \bv_{i,h}; \bmu) &=& - \mathcal{A}_i^{\pgd} (\mathbf{1}_{\Omega_i,h}^{k,j}, \bv_{i,h} ; \bmu)\,,\\
        \mathcal{B}_i^{\pgd}(q_{i,h},\bu^{k,j}_{i,\pgd}(\bmu); \bmu) + \mathcal{C}_i (p^{k,j}_{i,\pgd}(\bmu),q_{i,h}) &=& - \mathcal{B}_i^{\pgd} (q_{i,h},\mathbf{1}_{\Omega_i,h}^{k,j}; \bmu)\,,
	\end{eqnarray}
\end{subequations}
for all $(\bv_{i,h}, q_{i,h}) \in \bV_i^h \times W_i^h$ and for all $\bmu \in \mathcal{P}$.

\smallskip

We remark that the functions
\begin{equation}\label{eq:stokesPGDbasis} 
\{\bu^{k,j}_{i,\pgd}(\bmu),\, p^{k,j}_{i,\pgd}(\bmu)\}_{j = 1, \dots, N_{\Gamma_i}, \, k = 1, \dots, d}
\end{equation}
computed in \eqref{eq:StokesLocalPGDUnity} form a \textit{PGD expansion} for the Stokes problem in $\Omega_i$.  
Together with ${\bu}^{0,f}_{i,\pgd}(\bmu)$ and ${p}^f_{i,\pgd}(\bmu)$, the PGD expansion permits to obtain the local surrogate models for the local problem \eqref{eq:StokesMultiDomainGalerkin} through computationally cheap linear combinations analogous to \eqref{eq:upPGDLambda}. We detail this procedure in Sect.~\ref{sect:StokesAScoupling}.

\subsection{Online surrogate-based coupling procedure for the Stokes problem}
\label{sect:StokesAScoupling}

In the online phase, the local surrogate models computed in the offline phase (see Sect.~\ref{sect:StokesSurrogateModel}) are coupled to obtain the global solution of the Stokes problem \eqref{eq:StokesGlobal} for any given set of parameters $\bar{\bmu} \in \mathcal{P}$. The coupling is performed by solving the interface system \eqref{eq:StokesMultiDomainAlgebraicInterface} after replacing the solution of the local algebraic problems \eqref{eq:localStokesAlgebraicData} and \eqref{eq:localStokesAlgebraicLambda}  by the evaluation of suitable surrogate models for the assigned set of parameters $\bar{\bmu} \in \mathcal{P}$.

More precisely, instead of solving the linear systems \eqref{eq:localStokesAlgebraicData}, we evaluate the local surrogate models computed in \eqref{eq:StokesLocalPGDData} at $\bar{\bmu}$ to obtain $\bu^{0,f}_{i,\pgd} (\bar{\bmu})$ and $p^f_{i,\pgd} (\bar{\bmu})$, for $i=1,2$. Moreover, for $i=1,2$, let us introduce the PGD-based Stokes operators
\begin{subequations}\label{eq:stokesPGDoperators_up}
\begin{eqnarray}
&& \mathcal{S}_{u,i}^\pgd : \bLambda_{\Gamma_i} \to \sum_{k=1}^d \sum_{j=1}^{N_{\Gamma_i}} \Lambda_i^{k,j}(\bar{\bmu}) \, \bu^{k,j}_{i,\pgd}(\bar{\bmu}) + \mat{E}_{\Gamma_i \to \Omega_i} \bLambda_{\Gamma_i}, \label{eq:stokesPGDoperators_u} \\
&& \mathcal{S}_{p,i}^\pgd : \bLambda_{\Gamma_i} \to \sum_{k=1}^d \sum_{j=1}^{N_{\Gamma_i}} \Lambda_i^{k,j}(\bar{\bmu}) \, p^{k,j}_{i,\pgd}(\bar{\bmu}), \label{eq:stokesPGDoperators_p}
\end{eqnarray}
\end{subequations}
where $\bu^{k,j}_{i,\pgd}(\bar{\bmu})$ and $p^{k,j}_{i,\pgd}(\bar{\bmu})$ are the local PGD expansion functions \eqref{eq:stokesPGDbasis} at $\bar{\bmu}$. These operators perform linear combinations of the precomputed local surrogate models for any given vector of coefficient $\bLambda_{\Gamma_i}$, and they are used to replace the solution of problems \eqref{eq:localStokesAlgebraicLambda}.

Therefore, the interface system \eqref{eq:StokesMultiDomainAlgebraicInterface} can be reformulated as: for any assigned set of parameters $\bar{\bmu} \in \mathcal{P}$, find $\bLambda_{\Gamma_1}$ and $\bLambda_{\Gamma_2}$ such that
\begin{equation}\label{eq:stokesInterfacePGD}
 \begin{pmatrix}
  \mat{I}_{\Gamma_1\Gamma_1} & - \mat{R}_{\Omega_2 \to \Gamma_1} \, \mathcal{S}_{u,2}^\pgd \\
  - \mat{R}_{\Omega_1 \to \Gamma_2} \, \mathcal{S}_{u,1}^\pgd & \mat{I}_{\Gamma_2\Gamma_2}
 \end{pmatrix}
 \begin{pmatrix}
  \boldsymbol{\Lambda}_{\Gamma_1} \\
  \boldsymbol{\Lambda}_{\Gamma_2}
 \end{pmatrix} \\
= 
 \begin{pmatrix}
  \mathbf{g}_{\Gamma_1} \\
  \mathbf{g}_{\Gamma_2}
 \end{pmatrix}
+
\begin{pmatrix}
    \mat{R}_{\Omega_2 \to \Gamma_1} \, \bu^{0,f}_{2,\pgd} (\bar{\bmu}) \\
    \mat{R}_{\Omega_1 \to \Gamma_2} \, \bu^{0,f}_{1,\pgd} (\bar{\bmu})
\end{pmatrix}.
\end{equation}

System \eqref{eq:stokesInterfacePGD} is still solved iteratively using GMRES, and, at convergence, the values ${\bLambda}_{\Gamma_1}$ and $\bLambda_{\Gamma_2}$ ensure the smooth coupling of the local PGD solutions, and the global PGD approximation of the Stokes problem \eqref{eq:StokesGlobal} is defined as 
\begin{subequations}\label{eq:pgd_velPressure}
\begin{eqnarray}
    \bu^\texttt{PGD}(\bar{\bmu}) &=& 
    \begin{cases}
        \mathcal{S}_{u,1}^\pgd \, \bLambda_{\Gamma_1} + \bu^{0,f}_{1,\pgd}(\bar{\bmu}) + 
        \bg_{\Omega_1}^D(\bar{\bmu}) \quad &\text{in } \Omega_1 ,\\[5pt]
        \mathcal{S}_{u,2}^\pgd \, \bLambda_{\Gamma_2} + \bu^{0,f}_{2,\pgd}(\bar{\bmu}) + 
        \bg_{\Omega_2}^D(\bar{\bmu}) \quad &\text{in } \Omega_2 \setminus \Omega_{12}, \\
    \end{cases} \\
    p^\texttt{PGD}(\bar{\bmu}) &=& 
    \begin{cases}
        \mathcal{S}_{p,1}^\pgd \, \bLambda_{\Gamma_1} + p^f_{1,\pgd}(\bar{\bmu}) \quad &\text{in } \Omega_1, \\[5pt]
        \mathcal{S}_{p,2}^\pgd \, \bLambda_{\Gamma_2} + p^f_{2,\pgd}(\bar{\bmu}) \quad &\text{in } \Omega_2 \setminus \Omega_{12}.
    \end{cases}
\end{eqnarray}
\end{subequations}

The overall computational cost of the online phase is due to two main components:
\begin{enumerate}
    \item The evaluation of the PGD expansion \eqref{eq:stokesPGDbasis} for the selected parameter $\bar{\bmu}$. This involves computing the sum of the spatial and parametric modes \eqref{eq:SepSolnVelPress}, and this is done before starting the iterations to solve system \eqref{eq:stokesInterfacePGD}.
    \item The computation of the PGD-based operators $\mathcal{S}_{u,i}^\pgd$. This corresponds to performing the linear combinations \eqref{eq:stokesPGDoperators_u} for any given set of coefficients $\bLambda_{\Gamma_i}$ at each iteration of the GMRES method.
\end{enumerate}
As no solution of local problems is involved, the cost of the DD-PGD approach is much cheaper than working on \eqref{eq:StokesMultiDomainAlgebraicInterface} directly.
    
\section{The parametric Stokes-Darcy problem}
\label{sect:Darcy}

This section is devoted to the Stokes-Darcy coupled parametric problem. To avoid introducing further notation, the same letters used in Sect.~\ref{sect:Stokes} will still denote velocity, pressure and external forces/boundary data. However, throughout this section, all quantities labelled with a lower index `1' are associated with the Stokes problem, while the lower index `2' denotes physical quantities in the Darcy domain.

\subsection{Two-domain formulation of the parametric Stokes-Darcy problem}\label{sec:twoDomainStokesDarcy}

To couple the Stokes and the Darcy models, we consider the ICDD concept \cite{Discacciati:2024:JCP} that uses a thin overlapping region between the free fluid and the porous medium to model the transition between the two fluid regimes. In this framework, the problem is formulated as follows. Let $\Omega_1, \Omega_2 \subset \mathbb{R}^d$ ($d=2,3$) be the Stokes and the Darcy domains, respectively, with $\Omega=\Omega_1 \cup \Omega_2$, and $\Omega_{12} = \Omega_1 \cap \Omega_2 \not= \emptyset$ being the transition region. We use the same notation as in Sect.~\ref{sec:StokesProbStatement} for the parts of the boundary where different boundary conditions are imposed and for the interfaces,  and we still assume that the interfaces do not intersect, i.e., $dist(\Gamma_1,\Gamma_2)>0$. The two-domain formulation of the parametric Stokes-Darcy problem then becomes: for all $\bmu \in \mathcal{P}$, find the Stokes velocity $\bu_1(\bmu)$ and pressure $p_1(\bmu)$, and the Darcy velocity $\bu_2(\bmu)$ and pressure $p_2(\bmu)$ such that
\begin{subequations}
	\label{eq:StokesD1}
	\begin{eqnarray}
		-\nabla \cdot \bsigma(\bu_1(\bmu), p_1(\bmu); \bmu) = \beff_1(\bmu) &&\quad \text{in } \Omega_1,\\
		\nabla \cdot \bu_1(\bmu) = 0 &&\quad \text{in } \Omega_1,\\
		\bu_1(\bmu) = \bg^D_1(\bmu) &&\quad \text{on } \Gamma^D_1,\\
		\bsigma(\bu_1(\bmu), p_1(\bmu); \bmu) \bn = \bg^N_1(\bmu) &&\quad \text{on } \Gamma^N_1,
	\end{eqnarray}
\end{subequations}
and
\begin{subequations}
	\label{eq:DarcyD2}
	\begin{eqnarray}
		\nu(\bmu)\mat{K}(\bmu)^{-1}\bu_2(\bmu) + \nabla p_2(\bmu) = \beff_2(\bmu) &&\quad \text{in } \Omega_2,\\
		\nabla \cdot \bu_2(\bmu) = 0 &&\quad \text{in } \Omega_2,\\
		p_2(\bmu) = g_2^D(\bmu) &&\quad \text{on } \Gamma_2^D,\\
		\bu_2(\bmu) \cdot \bn = g_2^N(\bmu) &&\quad \text{on } \Gamma_2^N,
	\end{eqnarray}
\end{subequations}
with the coupling conditions
\begin{subequations}\label{eq:couplingConditionSD}
\begin{eqnarray}
    \bu_1(\bmu) = \bu_2(\bmu) && \quad \text{on } \Gamma_1, \label{eq:couplingConditionSD1}\\
    p_1(\bmu) = p_2(\bmu) && \quad \text{on } \Gamma_2. \label{eq:couplingConditionSD2}
\end{eqnarray}
\end{subequations}
The regular enough functions $\beff_1(\bmu),\ \bg^D_1(\bmu),\ \bg^N_1(\bmu),\ \beff_2(\bmu),\ g_2^D(\bmu)$ and $g_2^N(\bmu)$ represent external forces and boundary conditions, and $\mat{K}(\bmu)$ is the parametric permeability tensor such that, for all $\bmu \in \mathcal{P}$, $\bx^T\mat{K}(\bmu)\bx > 0$ for all $\bx \in \mathbb{R}^d$, $\bx \not= \mathbf{0}$.
The global velocity and pressure are defined as in \eqref{eq:globalVelocityPressure}.

Using \eqref{eq:couplingConditionSD1} as a boundary condition for the Stokes problem \eqref{eq:StokesD1}, the latter can be treated following the same procedure described in Sect.~\ref{sect:Stokes}. Therefore, we focus here on Darcy's problem and we consider the discrete case using computational meshes that coincide in the transition region $\Omega_{12}$ and are conforming with the interfaces like in the Stokes case (see Sect.~\ref{sec:StokesProbStatement}). In general, the Stokes and the Darcy problems are discretized using different finite element spaces, e.g., Taylor-Hood for Stokes and Raviart-Thomas elements for Darcy \cite{Boffi:2013}. Our approach can easily accommodate this situation because the local Stokes and Darcy problems exchange data only through the interface conditions  \eqref{eq:couplingConditionSD}, so that interpolation operators can be easily used to transfer information between one interface grid and another. However, for the sake of simplicity of exposition, we assume that the same continuous linear finite element basis functions are used to discretize both problems. For this purpose, we consider the stabilized formulation for Stokes introduced in Sect.~\ref{sec:StokesProbStatement}, while for Darcy we adopt the Masud-Hughes stabilized form \cite{Masud:2002:CMAME}. More precisely, let
\begin{equation*}
\begin{array}{c}
    \widetilde{\bV}_2^h = \{ \bv \in [X_2^1]^d \; : \; \bv \cdot \bn = 0 \text{ on } \Gamma_2^N \cup \Gamma_2 \},
    \quad
    \widetilde{W}_2^h = X_2^1 \cap L^2(\Omega_2),
    \\
    \widetilde{Y}_2^h = \{ \lambda \in C^0(\Gamma_2) \, : \, \lambda = 0 \text{ at } \overline{\Gamma}_2^D \cap \overline{\Gamma}_2 \text{ if } \overline{\Gamma}_2^D \cap \overline{\Gamma}_2 \not= \emptyset, \text{ and } \exists \, w \in \widetilde{W}_2^h \text{ s.t } w=\lambda \text{ on } \Gamma_2 \}.
\end{array}
\end{equation*}
Then, at the discrete level, we can decompose the Darcy problem \eqref{eq:DarcyD2} by introducing the following problems (using the strong form for the sake of clarity):
\begin{enumerate}
    \item Darcy problem depending on assigned data: find $\bu_{2,h}^f(\bmu)$ and $p_{2,h}^f(\bmu)$ such that
    \begin{subequations}
	\label{eq:DarcyD2Data}
	\begin{eqnarray}
		\nu(\bmu)\mat{K}(\bmu)^{-1}\bu_{2,h}^f(\bmu) + \nabla p_{2,h}^f(\bmu) = \beff_2(\bmu) &&\quad \text{in } \Omega_2,\\
		\nabla \cdot \bu_{2,h}^f(\bmu) = 0 &&\quad \text{in } \Omega_2,\\
		p_{2,h}^f(\bmu) = g_2^D(\bmu) &&\quad \text{on } \Gamma_2^D,\\
		\bu_{2,h}^f(\bmu) \cdot \bn = g_2^N(\bmu) &&\quad \text{on } \Gamma_2^N, \\
		p_{2,h}^f(\bmu) = \widetilde{g}_2^D(\bmu) &&\quad \text{on } \Gamma_2\be{,}
	\end{eqnarray}
\end{subequations}
where $\widetilde{g}_2^D(\bmu)$ on $\Gamma_2$ is introduced to avoid possible discontinuities in the pressure near the interface $\Gamma_2$. This auxiliary function is defined as
\begin{equation*}
    \widetilde{g}_2^D(\bmu) =
    \left\{
    \begin{array}{ll}
    \widehat{g}_2^D(\bmu) & \text{ if } \overline{\Gamma}_2^D \cap \overline{\Gamma}_2 \not= \emptyset , \\
    0 & \text{ if } \overline{\Gamma}_2^D \cap \overline{\Gamma}_2 = \emptyset,
    \end{array}
    \right.
\end{equation*}
where $\widehat{g}_2^D(\bmu)$ is a suitable continuous prolongation of $g_2^D(\bmu)$ on $\Gamma_2$.
Letting $\bg_{\Omega_2}^N (\bmu) \in [X_2^1]^d$ be a suitable continuous extension of $g_2^N(\bmu)$ such that $\bg_{\Omega_2}^N (\bmu) \cdot \bn = g_2^N(\bmu)$ at the dofs on $\Gamma_2^N$, we can split
\begin{equation*}
    \bu_{2,h}^f(\bmu) = \bu_{2,h}^{0,f}(\bmu) + \bg_{\Omega_2}^N (\bmu).
\end{equation*}

   \item Darcy problem depending on auxiliary interface data: find $\bu_{2,h}^\lambda (\bmu)$ and $p_{2,h}^\lambda (\bmu)$ such that
   \begin{subequations}
	\label{eq:DarcyD2Lambda}
	\begin{eqnarray}
		\nu(\bmu)\mat{K}(\bmu)^{-1}\bu_{2,h}^\lambda(\bmu) + \nabla p_{2,h}^\lambda(\bmu) = \mathbf{0} &&\quad \text{in } \Omega_2,\\
		\nabla \cdot \bu_{2,h}^\lambda(\bmu) = 0 &&\quad \text{in } \Omega_2,\\
		p_{2,h}^\lambda(\bmu) = 0 &&\quad \text{on } \Gamma_2^D,\\
		\bu_{2,h}^\lambda(\bmu) \cdot \bn = 0 &&\quad \text{on } \Gamma_2^N, \\
		p_{2,h}^\lambda(\bmu) = \lambda_{2,h}(\bmu) && \quad \text{on } \Gamma_2,
	\end{eqnarray}
\end{subequations}
where the auxiliary function $\lambda_{2,h}(\bmu) \in \widetilde{Y}_2^h$ is chosen in such a way that
\begin{equation*}
    \bu_{2,h}(\bmu) = \bu_{2,h}^f (\bmu) + \bu_{2,h}^\lambda(\bmu)
    \quad \text{and} \quad
    p_{2,h}(\bmu) = p_{2,h}^f (\bmu) + p_{2,h}^\lambda (\bmu)
    \quad \text{in } \Omega_2,
\end{equation*}
with $\bu_{2,h}(\bmu)$ and $p_{2,h}(\bmu)$ being the Galerkin approximations of the Darcy velocity and pressure in \eqref{eq:DarcyD2}.
\end{enumerate}

The discrete counterpart of the coupling conditions \eqref{eq:couplingConditionSD} can be written as
\begin{subequations}\label{eq:couplingConditionSDDiscrete}
    \begin{eqnarray}
        \blambda_{1,h}(\bmu) - \bu_{2,h}^0 (\bmu)_{\vert\Gamma_1} = (-\bg_{\Omega_1}^D (\bmu) + \bg_{\Omega_2}^N (\bmu))_{\vert\Gamma_1} && \text{on } \Gamma_1, \\
        \lambda_{2,h}(\bmu) - p_{1,h}(\bmu)_{\vert\Gamma_2} = -\widetilde{g}_2^D(\bmu) && \text{on } \Gamma_2,
    \end{eqnarray}
\end{subequations}
where we denoted
\begin{equation*}
    \bu_{2,h}^0 (\bmu) = \bu_{2,h}^{0,f}(\bmu) + \bu_{2,h}^\lambda (\bmu)\,.
\end{equation*}
Following \cite{Masud:2002:CMAME}, we define the bilinear forms and functionals: for all $\bu,\bv \in \widetilde{\bV}_2^h$, $q \in \widetilde{W}_2^h$, and $\bmu \in \mathcal{P}$,
\begin{eqnarray*}
    \widetilde{\mathcal{A}}_2 (\bu,\bv;\bmu) &=& \frac{1}{2} \int_{\Omega_2} \left( \nu (\bmu) \mat{K}^{-1} (\bmu) \bu \right) \cdot \bv \, d\bx \\
    \widetilde{\mathcal{B}}_2 (q, \bv) &=& - \int_{\Omega_2} q \, \nabla \cdot \bv \, d\bx - \frac{1}{2} \int_{\Omega_2} \nabla q \cdot \bv \, d\bx , \\
    \widetilde{\mathcal{C}}_2 (p,q;\bmu) &=& - \frac{1}{2} \int_{\Omega_2} \nabla q \cdot (\nu^{-1}(\bmu) \mat{K}(\bmu) \nabla p) \, d\bx, \\
    \widetilde{\mathcal{F}}_2 (\bv ; \bmu) &=& \frac{1}{2} \int_{\Omega_2} \beff_2(\bmu) \cdot \bv \, d\bx - \int_{\Gamma_2^D} g_2^{D}(\bmu) (\bv\cdot\bn)\, d\bx, \\
    \widetilde{\mathcal{G}}_2 (q;\bmu) &=& - \frac{1}{2} \int_{\Omega_2} \nabla q \cdot (\nu^{-1}(\bmu) \mat{K}(\bmu) \beff_2(\bmu) ) \, d\bx.
\end{eqnarray*}
Then, the Galerkin approximation of the Stokes-Darcy problem \eqref{eq:StokesD1}-\eqref{eq:DarcyD2} with the coupling conditions expressed as in \eqref{eq:couplingConditionSDDiscrete} becomes: for all $\bmu \in \mathcal{P}$, find $(\bu_{1,h}^0(\bmu),p_{1,h}(\bmu)) \in \bV_1^h \times W_1^h$, $(\bu_{2,h}^0(\bmu),p_{2,h}(\bmu)) \in \widetilde{\bV}_2^h \times \widetilde{W}_2^h$, $\blambda_{1,h}(\bmu) \in \bY_1^h$ and $\lambda_{2,h}(\bmu) \in \widetilde{Y}_2^h$ such that, for all $(\bv_{1,h},q_{1,h}) \in \bV_1^h \times W_1^h$ and $(\bv_{2,h},q_{2,h}) \in \widetilde{\bV}_2^h \times \widetilde{W}_2^h$,
\begin{subequations}\label{eq:StokesDarcyTwoDomainGalerkin}
	\begin{eqnarray}
		\hspace*{-5mm}
		\mathcal{A}_1(\bu_{1,h}^{0}(\bmu) + \blambda_{\Omega_1,h}(\bmu), \bv_{1,h}; \bmu) &+& \mathcal{B}_1(p_{1,h}(\bmu), \bv_{1,h}; \bmu) \nonumber \\ 
        \hspace*{-5mm}
        &=& \mathcal{F}_1(\bv_{1,h};\bmu) - \mathcal{A}_1 (\bg_{\Omega_1}^D(\bmu), \bv_{1,h} ; \bmu),\\[5pt]
        \hspace*{-5mm}
        \mathcal{B}_1(q_{1,h},\bu_{1,h}^0(\bmu) + \blambda_{\Omega_1,h}(\bmu); \bmu) + \mathcal{C}_1 (p_{1,h}(\bmu),q_{1,h}) &=& \mathcal{G}_1 (q_{1,h};\bmu) - \mathcal{B}_1 (q_{1,h},\bg_{\Omega_1}^D(\bmu); \bmu)\,, \\[5pt]
        \hspace*{-5mm}
        \widetilde{\mathcal{A}}_2(\bu_{2,h}^{0}(\bmu), \bv_{2,h}; \bmu) + \widetilde{\mathcal{B}}_2(p_{2,h}(\bmu), \bv_{2,h}) &+& \int_{\Gamma_2} \lambda_{2,h}(\bmu) (\bv_{2,h}\cdot \bn) \, d\bx \nonumber \\
        \hspace*{-5mm}
        &=& \widetilde{\mathcal{F}}_2(\bv_{2,h};\bmu) - \widetilde{\mathcal{A}}_2 (\bg_{\Omega_2}^N(\bmu), \bv_{2,h} ; \bmu),\\[5pt]
        \hspace*{-5mm}
        \widetilde{\mathcal{B}}_2(q_{2,h},\bu_{2,h}^0(\bmu)) + \widetilde{\mathcal{C}}_2 (p_{2,h}(\bmu),q_{2,h};\bmu) &=& \widetilde{\mathcal{G}}_2 (q_{2,h};\bmu) - \widetilde{\mathcal{B}}_2 (q_{2,h},\bg_{\Omega_2}^N(\bmu))\,, \\[5pt]
        \hspace*{-5mm}
        \blambda_{1,h}(\bmu) - \bu^0_{2,h}(\bmu)_{\vert\Gamma_1} &=& (- \bg^D_{\Omega_1}(\bmu) + \bg^N_{\Omega_2}(\bmu))_{\vert\Gamma_1} \quad \text{on } \Gamma_1, \\[5pt]
        \hspace*{-5mm}
        \lambda_{2,h}(\bmu) - p_{1,h}(\bmu)_{\vert\Gamma_2} &=& -\widetilde{g}_2^D(\bmu) \quad \text{on } \Gamma_2.
	\end{eqnarray}
\end{subequations}

Using an analogous notation to the Stokes-Stokes coupling (see \eqref{eq:StokesMultiDomainAlgebraic}), the algebraic form of the two-domain Stokes-Darcy problem \eqref{eq:StokesDarcyTwoDomainGalerkin} can be written as
\begin{equation}\label{eq:StokesDarcyTwoDomainAlgebraic}
 \begin{pmatrix}
 \mat{A}^1_{II} & (\mat{B}^1_{\all I})^T & \mat{0} & \mat{0} & \mat{A}^1_{I\Gamma_1} & \mat{0} \\[2pt]
 \mat{B}^1_{\all I} & \mat{C}^1_{\all \all} & \mat{0} & \mat{0} & \mat{B}^1_{\all \Gamma_1} & \mat{0} \\[2pt]
 \mat{0} & \mat{0} & \widetilde{\mat{A}}^2_{II} & (\widetilde{\mat{B}}^2_{\all I})^T & \mat{0} & \mat{M}^2_{I\Gamma_2} \\[2pt]
 \mat{0} & \mat{0} & \widetilde{\mat{B}}^2_{\all I} & \widetilde{\mat{C}}^2_{\all \all} & \mat{0} & \mat{0} \\[2pt]
 \mat{0} & \mat{0} & -\mat{R}_{\Omega_2 \to \Gamma_1} & \mat{0} & \mat{I}_{\Gamma_1\Gamma_1} & \mat{0} \\[2pt]
 \mat{0} & -\mat{R}_{\Omega_1\to\Gamma_2} & \mat{0} & \mat{0} & \mat{0} & \mat{I}_{\Gamma_2\Gamma_2}
 \end{pmatrix}
 \begin{pmatrix}
     \mathbf{u}^1_I \\[2pt]
     \mathbf{p}^1_\all \\[2pt]
     \widetilde{\mathbf{u}}^2_I \\[2pt]
     \widetilde{\mathbf{p}}^2_\all \\[2pt]
     \boldsymbol{\Lambda}_{\Gamma_1} \\[2pt]
     \widetilde{\boldsymbol{\Lambda}}_{\Gamma_2}
 \end{pmatrix}
 =
 \begin{pmatrix}
     \mathbf{f}^1_I \\[2pt]
     \mathbf{g}^1_\all \\[2pt]
     \widetilde{\mathbf{f}}^2_I \\[2pt]
     \widetilde{\mathbf{g}}^2_\all \\[2pt]
     \widetilde{\mathbf{g}}_{\Gamma_1} \\[2pt]
     \widetilde{\mathbf{g}}_{\Gamma_2}
 \end{pmatrix}.
\end{equation}

In this case, the equivalent interface system becomes
\begin{equation}\label{eq:StokesDarcyAlgebraicInterface}
\begin{array}{rcl}
\left(
 \begin{pmatrix}
  \mat{I}_{\Gamma_1\Gamma_1} & \mat{0} \\
  \mat{0} & \mat{I}_{\Gamma_2\Gamma_2}
 \end{pmatrix}
\right. 
&-&
 \begin{pmatrix}
  \mat{0} & \mat{0} \\
  \mat{0} & \mat{R}_{\Omega_1\to\Gamma_2}
 \end{pmatrix}
 \begin{pmatrix}
  \mat{A}^1_{II} & (\mat{B}^1_{\all I})^T \\
  \mat{B}^1_{\all I} & \mat{C}^1_{\all \all} 
 \end{pmatrix}^{-1}
 \begin{pmatrix}
  -\mat{A}^1_{I\Gamma_1} & \mat{0} \\
  -\mat{B}^1_{\all \Gamma_1} & \mat{0} 
 \end{pmatrix} \\
&-&
\left.
 \begin{pmatrix}
  \mat{R}_{\Omega_2 \to \Gamma_1} & \mat{0} \\
  \mat{0} & \mat{0} 
 \end{pmatrix}
 \begin{pmatrix}
  \widetilde{\mat{A}}^2_{II} & (\widetilde{\mat{B}}^2_{\all I})^T \\
  \widetilde{\mat{B}}^2_{\all I} & \widetilde{\mat{C}}^2_{\all \all}
 \end{pmatrix}^{-1}
 \begin{pmatrix}
  \mat{0} & -\mat{M}^2_{I\Gamma_2} \\
  \mat{0} & \mat{0}
 \end{pmatrix}
\right)
\begin{pmatrix}
  \boldsymbol{\Lambda}_{\Gamma_1} \\
  \widetilde{\boldsymbol{\Lambda}}_{\Gamma_2}
\end{pmatrix} \\
= 
 \begin{pmatrix}
  \widetilde{\mathbf{g}}_{\Gamma_1} \\
  \widetilde{\mathbf{g}}_{\Gamma_2}
 \end{pmatrix}
&+&
 \begin{pmatrix}
  \mat{0} & \mat{0} \\
  \mat{0} & \mat{R}_{\Omega_1\to\Gamma_2}
 \end{pmatrix}
 \begin{pmatrix}
  \mat{A}^1_{II} & (\mat{B}^1_{\all I})^T \\
  \mat{B}^1_{\all I} & \mat{C}^1_{\all \all} 
 \end{pmatrix}^{-1}
 \begin{pmatrix}
  \mathbf{f}^1_I \\
  \mathbf{g}^1_\all \\
 \end{pmatrix} \\
&+&
 \begin{pmatrix}
  \mat{R}_{\Omega_2 \to \Gamma_1} & \mat{0} \\
  \mat{0} & \mat{0} 
 \end{pmatrix}
 \begin{pmatrix}
  \widetilde{\mat{A}}^2_{II} & (\widetilde{\mat{B}}^2_{\all I})^T \\
  \widetilde{\mat{B}}^2_{\all I} & \widetilde{\mat{C}}^2_{\all \all}
 \end{pmatrix}^{-1}
 \begin{pmatrix}
  \widetilde{\mathbf{f}}^2_I \\
  \widetilde{\mathbf{g}}^2_\all
\end{pmatrix}.
\end{array}
\end{equation}

The interface system involves solving the same Stokes problems \eqref{eq:localStokesAlgebraicData} and \eqref{eq:localStokesAlgebraicLambda} in $\Omega_1$, and the Darcy problems
\begin{equation}\label{eq:darcyAlgebraic}
\begin{pmatrix}
  \widetilde{\mat{A}}^2_{II} & (\widetilde{\mat{B}}^2_{\all I})^T \\
  \widetilde{\mat{B}}^2_{\all I} & \widetilde{\mat{C}}^2_{\all \all}
 \end{pmatrix}^{-1}
 \begin{pmatrix}
  \widetilde{\mathbf{f}}^2_I \\
  \widetilde{\mathbf{g}}^2_\all
\end{pmatrix}
\quad \text{and} \quad
\begin{pmatrix}
  \widetilde{\mat{A}}^2_{II} & (\widetilde{\mat{B}}^2_{\all I})^T \\
  \widetilde{\mat{B}}^2_{\all I} & \widetilde{\mat{C}}^2_{\all \all}
 \end{pmatrix}^{-1}
 \begin{pmatrix}
  -\mat{M}^2_{I\Gamma_2} \\
  \mat{0}
 \end{pmatrix}
 \widetilde{\bLambda}_{\Gamma_2}
\end{equation}
that are, respectively, the counterparts of \eqref{eq:DarcyD2Data} and \eqref{eq:DarcyD2Lambda}. In the next section, we explain how to construct surrogate models for these problems using PGD in the offline phase of the algorithm.

\subsection{Local surrogate models using proper generalized decomposition for the Darcy problem}
\label{sect:DarcySurrogateModel}

Analogously to Sect.~\ref{sect:StokesSurrogateModel}, we assume that the data of the Darcy problem is expressed in separable form. In particular, we write
\begin{equation*}
    \begin{array}{rclrcl}
    \nu(\bmu)^{-1} &=& \displaystyle\sum_{\ell=1}^{\widehat{N}_{\nu}} \widehat{\xi}_{\nu}^\ell(\bmu)\widehat{b}_{\nu}^\ell(\bx) \,,
    & \quad
    \beff_2(\bmu) &=& \displaystyle\sum_{\ell=1}^{N_{2,f}} \xi_{2,f}^\ell(\bmu) \bb_{2,f}^\ell(\bx) \,,
    \\
    \mat{K}(\bmu) &=& \displaystyle\sum_{\ell=1}^{N_{\kappa}} \xi_{\kappa}^\ell(\bmu) \bb_{\kappa}^\ell(\bx) \,,
    & \quad
    g_{2}^D(\bmu) &=& \displaystyle\sum_{\ell=1}^{N_{2,D}}  \xi_{2,D}^\ell(\bmu) b_{2,D}^\ell(\bx) \,,
    \\
    \mat{K}(\bmu)^{-1} &=& \displaystyle\sum_{\ell=1}^{\widehat{N}_{\kappa}} \widehat{\xi}_{\kappa}^\ell(\bmu) \widehat{\bb}_{\kappa}^\ell(\bx) \,,
    & \quad
	\bg_{\Omega_2}^N(\bmu) &=& \displaystyle\sum_{\ell=1}^{N_{2,N}}  \xi_{2,N}^\ell(\bmu) \bb_{2,N}^\ell(\bx) \,.
    \end{array}
\end{equation*}
\begin{rem}
    Note that if $\nu(\bmu)$ and $\mat{K}(\bmu)$ are given in separable form, $\nu(\bmu)^{-1}$ and $\mat{K}(\bmu)^{-1}$ will not typically be separable, at least analytically. One can circumvent this issue by numerically constructing a separable approximation of these quantities (see \cite{Diez:2020:ACME} and \ref{app:StokesDarcyData} for details).
\end{rem}
We focus now on problem \eqref{eq:DarcyD2Data}. We still consider approximations of the velocity $\bu_{2,h}^{0,f}$ and pressure $p_{2,h}^f$ in the form \eqref{eq:SepSolnVelPress} with the Dirichlet datum $\bg_{\Omega_2}^N(\bmu)$ accounted for by ad-hoc smooth modes. Assuming affine parameter dependence, for all $\bu,\bv \in \widetilde{\bV}_2^h$, $p,q \in \widetilde{W}_2^h$ and $\bmu \in \mathcal{P}$, we introduce the bilinear forms
\begin{eqnarray*}
    \widetilde{\mathcal{A}}_2^\pgd (\bu,\bv;\bmu) &=& \frac{1}{2} \sum_{\ell=1}^{N_\nu} \sum_{m=1}^{\widehat{N}_\kappa} \xi_\nu^\ell (\bmu) \widehat{\xi}_\kappa^m (\bmu) \int_{\Omega_2} (b_\nu^\ell (\bx) \widehat{\bb}_\kappa^m (\bx) \bu) \cdot \bv \, d\bx \,, \\
    \widetilde{\mathcal{C}}_2^\pgd (p,q;\bmu) &=& - \frac{1}{2} \sum_{\ell = 1}^{\widehat{N}_\nu} \sum_{m=1}^{N_\kappa} \widehat{\xi}_\nu^\ell (\bmu) \xi_\kappa^m (\bmu) \int_{\Omega_2} \nabla q \cdot (\widehat{b}_\nu^\ell (\bx) \bb_\kappa^m (\bx) \nabla p) \, d\bx \,,
\end{eqnarray*}
and the functionals
\begin{eqnarray*}
    \widetilde{\mathcal{F}}_2^\pgd (\bv;\bmu) &=& \frac{1}{2} \sum_{\ell=1}^{N_{2,f}} \xi_{2,f}^\ell (\bmu) \int_{\Omega_2} \bb_{2,f}^\ell (\bx) \cdot \bv \, d\bx - \sum_{\ell=1}^{N_{2,D}} \xi_{2,D}^\ell (\bmu) \int_{\Gamma_N} b_{2,D}^\ell (\bx) \, (\bv \cdot \bn) \, d\bx \, , \\
    \widetilde{\mathcal{G}}_2^\pgd (q;\bmu) &=& - \frac{1}{2} \sum_{\ell = 1}^{\widehat{N}_\nu} \sum_{m=1}^{N_\kappa} \sum_{n=1}^{N_{2,f}} \widehat{\xi}_\nu^\ell (\bmu) \xi_\kappa^m (\bmu) \xi_{2,f}^n (\bmu) \int_{\Omega_2} \nabla q \cdot \left(\widehat{b}_\nu^\ell (\bx) \bb_\kappa^m (\bx) \bb_{2,f}^n (\bx) \right) \, d\bx\,.
\end{eqnarray*}
Then, the PGD approximations $\bu_{2,\pgd}^{0,f}(\bmu)$ and $p_{2,\pgd}^f(\bmu)$ of the Darcy velocity $\bu_{2,h}^{0,f}(\bmu)$ and pressure $p_{2,h}^f(\bmu)$ in \eqref{eq:DarcyD2Data} are computed in the offline phase as the solution of the parametric problem:
\begin{subequations}\label{eq:problemPGDDarcyData}
    \begin{eqnarray}
     \widetilde{\mathcal{A}}_2^\pgd (\bu_{2,\pgd}^{0,f}(\bmu), \bv_{2,h};\bmu) &+& \widetilde{\mathcal{B}}_2 (p_{2,\pgd}^f(\bmu) , \bv_{2,h}) \nonumber \\
        &=& \widetilde{\mathcal{F}}_2^\pgd (\bv_{2,h}; \bmu)
        - \widetilde{\mathcal{A}}_2^\pgd ( \sum_{\ell=1}^{N_{2,N}} \xi_{2,N}^\ell (\bmu) \bb_{2,N}^\ell (\bx), \bv_{2,h} ; \bmu)\,, \\
        \widetilde{\mathcal{B}}_2 (q_{2,h}, \bu_{2,\pgd}^{0,f} (\bmu)) &+& \widetilde{\mathcal{C}}_2^\pgd (p_{2,\pgd}^f(\bmu),q_{2,h}; \bmu) \nonumber \\
        &=& \widetilde{\mathcal{G}}_2^\pgd (q_{2,h}; \bmu)
        - \widetilde{\mathcal{B}}_2 (q_{2,h}, \sum_{\ell = 1}^{N_{2,N}} \xi_{2,N}^\ell (\bmu) \bb_{2,N}^\ell (\bx)),
    \end{eqnarray}
\end{subequations}
for all $\bv_{2,h} \in \widetilde{\bV}_2^h$, $q_{2,h} \in \widetilde{W}_2^h$, and for all $\bmu \in \mathcal{P}$.

To deal with the generic trace function $\lambda_{2,h}(\bmu)$ on $\Gamma_2$ in \eqref{eq:DarcyD2Lambda}, we follow a similar approach to the one for the Stokes problem and express
\begin{equation*}
    \lambda_{2,h}(\bmu) = \lambda_{2,h}(\bx;\bmu) = \sum_{j=1}^{N_{\Gamma_2}} \widetilde{\Lambda}_2^j(\bmu) \, \eta_{2,h}^j(\bx)\,,
\end{equation*}
where now $\eta_{2,h}^j (\bx)$ are the basis functions of $\widetilde{Y}_2^h$ associated with the $N_{\Gamma_2}$ dofs on $\Gamma_2$, and $\widetilde{\Lambda}_2^j(\bmu)$ are real coefficients. Then, for all $\bmu \in \mathcal{P}$ and $j=1,\ldots,N_{\Gamma_2}$, we consider the auxiliary Darcy problems: find $(\bu_{2,h}^j (\bmu) , p_{2,h}^j (\bmu)) \in \widetilde{\bV}_2^h \times \widetilde{W}_2^h$ such that, for all $(\bv_{2,h},q_{2,h}) \in \widetilde{\bV}_2^h \times \widetilde{W}_2^h$,
\begin{eqnarray*}
    \widetilde{\mathcal{A}}_2 (\bu_{2,h}^j (\bmu) , \bv_{2,h} ; \bmu) + \widetilde{\mathcal{B}}_2 (p_{2,h}^j(\bmu), \bv_{2,h}) &=& - \int_{\Gamma_2} \eta_{2,h}^j \, (\bv_{2,h} \cdot \bn) \, d\bx \,, \\
    \widetilde{\mathcal{B}}_2 (q_{2,h}, \bu_{2,h}^j (\bmu)) + \widetilde{\mathcal{C}}_2 (p_{2,h}^j(\bmu), q_{2,h}; \bmu) &=& 0 \, .
\end{eqnarray*}
Owing to the linearity of the Darcy problem, by construction it holds
\begin{equation}\label{eq:darcyDecomposition}
    \bu_{2,h}^\lambda(\bmu) = \sum_{j=1}^{N_{\Gamma_2}} \widetilde{\Lambda}_2^j(\bmu) \, \bu_{2,h}^j (\bmu)
    \qquad \text{and} \qquad
    p_{2,h}^\lambda (\bmu) = \sum_{j=1}^{N_{\Gamma_2}} \widetilde{\Lambda}_2^j(\bmu) \, p_{2,h}^j (\bmu) \, .
\end{equation}
Then, the local surrogate models $\bu_{2,\pgd}^\lambda (\bmu)$ and $p_{2,\pgd}^\lambda (\bmu)$ can be obtained using linear combinations like \eqref{eq:darcyDecomposition}, but replacing the finite element velocity and pressure $\bu_{2,h}^j(\bmu)$ and $p_{2,h}^j(\bmu)$ by the \textit{PGD expansion}
\begin{equation}\label{eq:pdgBasisDarcy}
    \{ \bu_{2,\pgd}^j (\bmu), \, p_{2,\pgd}^j (\bmu) \}_{j=1,\ldots,N_{\Gamma_2}}.
\end{equation}
These are computed offline by solving the $N_{\Gamma_2}$ local PGD problems
\begin{subequations}\label{eq:problemDarcyBase}
\begin{eqnarray}
    \widetilde{\mathcal{A}}_2^\pgd (\bu_{2,\pgd}^j (\bmu) , \bv_{2,h} ; \bmu) + \widetilde{\mathcal{B}}_2 (p_{2,\pgd}^j(\bmu), \bv_{2,h}) &=& - \int_{\Gamma_2} \eta_{2,h}^j \, (\bv_{2,h} \cdot \bn) \, d\bx \,, \\
    \widetilde{\mathcal{B}}_2 (q_{2,h}, \bu_{2,\pgd}^j (\bmu)) + \widetilde{\mathcal{C}}_2^\pgd (p_{2,\pgd}^j(\bmu), q_{2,h}; \bmu) &=& 0 \, ,
\end{eqnarray}
\end{subequations}
for all $(\bv_{2,h},q_{2,h}) \in \widetilde{\bV}_2^h \times \widetilde{W}_2^h$, $\bmu \in \mathcal{P}$, and for $j=1,\ldots,N_{\Gamma_2}$. Note that these problems are the Darcy counterpart of the Stokes PGD problems \eqref{eq:StokesLocalPGDUnity}. In practice, problems \eqref{eq:problemPGDDarcyData} and \eqref{eq:problemDarcyBase} are solved using the non-intrusive Encapsulated PGD Algebraic Toolbox \cite{Diez:2020:ACME}.

\subsection{Online surrogate-based coupling for Stokes-Darcy}

In the online phase, we use the local surrogate models built in Sect.~\ref{sect:StokesSurrogateModel} for Stokes in $\Omega_1$ and in Sect.~\ref{sect:DarcySurrogateModel} for Darcy in $\Omega_2$ to obtain global surrogate Stokes-Darcy models for any given set of parameters $\bmu \in \mathcal{P}$. This will be achieved by iteratively solving the interface system \eqref{eq:StokesDarcyAlgebraicInterface} after replacing the solution of the local problems by the evaluation of precomputed local surrogate models. To this aim, for Stokes we will consider the same PGD operators defined in \eqref{eq:stokesPGDoperators_up}, while for Darcy we introduce the PGD operators
\begin{equation}\label{eq:darcyPGDoperators}
    \mathcal{D}_{u,2}^\pgd : \widetilde{\bLambda}_{\Gamma_2} \to \sum_{j=1}^{N_{\Gamma_2}} \widetilde{\Lambda}_2^{j} (\bar{\bmu}) \, \bu^{j}_{2,\pgd}(\bar{\bmu}),
    \qquad
    \mathcal{D}_{p,2}^\pgd : \widetilde{\bLambda}_{\Gamma_2} \to \sum_{j=1}^{N_{\Gamma_2}} \widetilde{\Lambda}_2^{j} (\bar{\bmu}) \, p^{j}_{2,\pgd}(\bar{\bmu}),
\end{equation}
that perform linear combinations with coefficients $\widetilde{\Lambda}_2^j(\bar{\bmu})$ of the PGD expansion functions \eqref{eq:pdgBasisDarcy}. Then, the interface system \eqref{eq:StokesDarcyAlgebraicInterface} can be expressed as: for any $\bar{\bmu} \in \mathcal{P}$, find $\bLambda_{\Gamma_1}$ and $\widetilde{\bLambda}_{\Gamma_2}$ such that
\begin{equation}\label{eq:stokesDarcyPGDinterface}
    \begin{pmatrix}
        \mat{I}_{\Gamma_1 \Gamma_1} & - \mat{R}_{\Omega_2 \to \Gamma_1} \mathcal{D}_{u,2}^\pgd \\
        -\mat{R}_{\Omega_1 \to \Gamma_2} \mathcal{S}_{p,1}^\pgd & \mat{I}_{\Gamma_2 \Gamma_2}
    \end{pmatrix}
    \begin{pmatrix}
        \bLambda_{\Gamma_1} \\
        \widetilde{\bLambda}_{\Gamma_2}
    \end{pmatrix}
    =
    \begin{pmatrix}
        \widetilde{\bg}_{\Gamma_1} \\
        \widetilde{\bg}_{\Gamma_2}
    \end{pmatrix}
    +
    \begin{pmatrix}
        \mat{R}_{\Omega_2\to\Gamma_1} \bu_{2,\pgd}^{0,f} (\bmu) \\
        \mat{R}_{\Omega_1 \to \Gamma_2} p_{1,\pgd}^f (\bmu)
    \end{pmatrix} \, .
\end{equation}
This linear system can be solved using GMRES and, at convergence, we can obtain the velocity and pressure of the global PGD surrogate model for Stokes-Darcy as
\begin{subequations}\label{eq:stokesDarcyPGDglobal}
\begin{eqnarray}
    \bu^\pgd (\bar{\bmu}) &=& 
    \begin{cases}
        \mathcal{S}_{u,1}^\pgd \, \bLambda_{\Gamma_1} + \bu^{0,f}_{1,\pgd}(\bar{\bmu}) + 
        \bg_{\Omega_1}^D(\bar{\bmu}) \quad &\text{in } \Omega_1 \\[5pt]
        \mathcal{D}_{u,2}^\pgd \, \widetilde{\bLambda}_{\Gamma_2} + \bu^{0,f}_{2,\pgd}(\bar{\bmu}) + 
        \bg_{\Omega_2}^N(\bar{\bmu}) \quad &\text{in } \Omega_2 \setminus \Omega_{12}, \\
    \end{cases} \\
    p^\texttt{PGD}(\bar{\bmu}) &=& 
    \begin{cases}
        \mathcal{S}_{p,1}^\pgd \, \bLambda_{\Gamma_1} + p^f_{1,\pgd}(\bar{\bmu}) \quad &\text{in } \Omega_1, \\[5pt]
        \mathcal{D}_{p,2}^\pgd \, \widetilde{\bLambda}_{\Gamma_2} + p^f_{2,\pgd}(\bar{\bmu}) \quad &\text{in } \Omega_2 \setminus \Omega_{12}.
    \end{cases}
\end{eqnarray}
\end{subequations}

Like in the Stokes case (see Sect. \ref{sect:StokesAScoupling}), the computational cost of the online phase is due to the evaluation of the Stokes and Darcy expansions \eqref{eq:stokesPGDbasis} and \eqref{eq:pdgBasisDarcy} at $\bar{\bmu}$, and to the computation of the linear combinations involved in the definition of the operators $\mathcal{S}_{p,1}^\pgd$ and $\mathcal{D}_{u,2}^\pgd$. The latter operation must be computed at each iteration of GMRES, while the former is performed once $\bar{\bmu}$ is selected, before starting the iterative method to solve system \eqref{eq:stokesDarcyPGDinterface}.

\section{Numerical results}
\label{sect:numericalResults}

In this section we present numerical tests\footnote{The numerical results presented in this section have been obtained using a PC with CPU Intel$^{\mbox{\tiny{\textregistered}}}$ Core\texttrademark\; i5-11400 @ 2.60GHz and 8GB RAM.} to assess the performance of the DD-PGD method in the context of multi-variable (i.e., velocity and pressure) and multi-physics (i.e., viscous flows and porous media) problems. We consider a Stokes-Stokes problem (Sect.~\ref{sect:StokesStokes}) and a Stokes-Darcy problem (Sect.~\ref{sec:Stokes-Darcy_synthetic}) both with analytic solution, and a Stokes-Darcy problem to simulate a cross-flow membrane filtration problem (Sect.~\ref{sect:StokesDarcyPhysical}) to test the method in the presence of physically relevant data.

To compute the local surrogate models, we use the Encapsulated PGD Algebraic Toolbox~\cite{Diez:2020:ACME} with the \texttt{mask} option to effectively handle multiple variables as explained in \ref{app:masks}. For all problems we employ tolerance $10^{-4}$ for the enrichment process and $10^{-3}$ for the PGD compression procedure. In the online phase, we use GMRES~\cite{Saad:1986:SISSC} without restart, with tolerance $10^{-6}$ on the relative residual to solve the system that is derived from the coupling conditions.

\subsection{Stokes-Stokes problem with analytic solution in two subdomains}\label{sect:StokesStokes}

Consider the Stokes problem \eqref{eq:StokesGlobal} in the spatial domain $\Omega = [0, 1] \times [0, 1]$, whose boundary is split into the disjoint subsets $\Gamma^N = [0, 1] \times \{0\}$ and $\Gamma^D = \partial\Omega \setminus \Gamma^N$. Let $\nu(\mu) = (1 - y) + y\mu$ be the space-dependent parametric viscosity, with $\mu \in [1,5]$ a scalar parameter. The force $\beff(\mu)$ and the boundary data $\bg^D(\mu)$ and $\bg^N(\mu)$ are chosen such that the exact solutions for velocity and pressure are
\begin{align*}
		\bu_\text{ex}(\mu) &= \frac{1}{100}\begin{bmatrix}
			3x - y\\
			3x^2 - 3y - x
		\end{bmatrix} + \begin{bmatrix}
		x^2(1 - x)^2(2y - 6y^2 + 4y^3)\\
		-y^2(1 - y)^2(2x - 6x^2 + 4x^3)
	\end{bmatrix}\mu, \\
		p_\text{ex}(\mu) &= y(3 - y) + x(1 - x^2)\mu.
	\end{align*}
We partition the domain $\Omega$ into two overlapping subdomains $\Omega_1 = [0, 0.55] \times [0, 1]$ and $\Omega_2 = [0.45, 1] \times [0, 1]$, so that the interfaces are $\Gamma_1 = \{ 0.55\} \times [0,1]$ and $\Gamma_2 = \{0.45\} \times [0,1]$. For the space discretization, we use the inf-sup stable $\mathbb{Q}_2-\mathbb{Q}_1$ Lagrangian finite elements \cite{Boffi:2013} on a uniform mesh with local mesh size $h = 5 \times 10^{-2}$, while the parametric domain $\mathcal{P}$ is discretized using pointwise collocation with spacing $h_\mu = 10^{-3}$.

Since the velocity is discretized using $\mathbb{Q}_2$ elements and considering the Dirichlet boundary condition on $\Gamma^D$, there are 40 dofs for each component of the velocity on each interface $\Gamma_i$ ($i=1,2$). Therefore, in each subdomain one has to solve 80 independent problems \eqref{eq:StokesLocalPGDUnity} and the problem \eqref{eq:StokesLocalPGDData}, which can be solved in parallel with a total offline computational cost of approximately $120$~s (including PGD compression). The total number of PGD modes computed after compression is 244 in $\Omega_1$ and 245 in $\Omega_2$.

In the online phase, taking, e.g., $\bar{\mu} =3$, the solution of the interface system \eqref{eq:stokesInterfacePGD} is computed in 27 GMRES iterations with a computational time of approximately $0.04$~s for the iterations only, and a total online phase time of $0.15$~s when factoring in pre-computations to setup the iterative system for the given value of the parameter $\bar{\mu}$ (see Sect. \ref{sect:StokesAScoupling}).
In comparison, a global FE solver computes the solution of this problem for a fixed value of the parameter $\mu$ in $0.19$~s, while solving the interface problem \eqref{eq:StokesMultiDomainAlgebraicInterface} using the classical DD-FEM approach also requires 27 GMRES iterations but with the computational time of $5.50$~s. Therefore, the DD-PGD approach has a computational cost comparable to a single classical FE solution, but it has a significant speed-up compared to the standard DD-FEM approach, with the computation time being 36 times faster, even when including the pre-computations. Moreover, the results of the offline phase can be further exploited for new parameter evaluations, without the need for new computations or re-training activities.

To quantify the accuracy of the method, we compute the relative error in the norm $L^2(\Omega)$ between the exact solution $v_{\text{ex}}$, the global finite element solution $v_\Omega^h$, the solution $v^h$ obtained by DD-FEM, and the DD-PGD solution $v^\pgd$ computed as in \eqref{eq:pgd_velPressure}, where $v$ denotes either the pressure or any component of the velocity. 
The numerical errors reported in Table \ref{tab:StokesErrors} show that the relative error for the velocity components computed by the DD-PGD procedure is slightly higher than its counterparts obtained by finite elements, while the pressure is approximated with similar accuracy.

\begin{figure}[bht]
    \begin{tabular}{ccc}
    \includegraphics[width=0.32\textwidth]{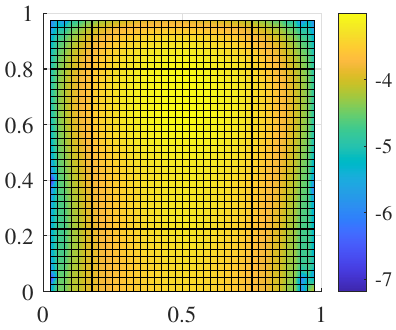} &
    \includegraphics[width=0.33\textwidth]{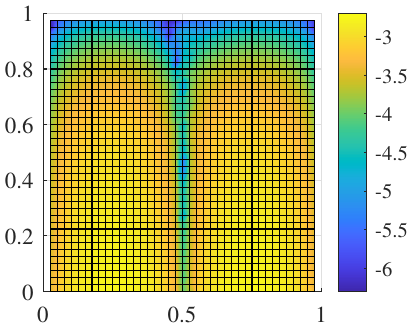} &
    \includegraphics[width=0.32\textwidth]{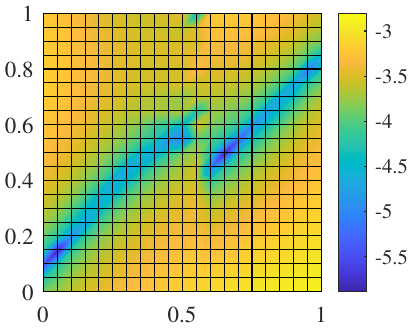} \\
    \footnotesize{Velocity ($x$ component)} &
    \footnotesize{Velocity ($y$ component)} &
    \footnotesize{Pressure}
    \end{tabular}
    \caption{For $\mu=3$, plot of $\log_{10}(|v^\pgd(\mu) - v_\text{ex}(\mu)| / \max_\Omega |v_\text{ex}(\mu)|)$ for $v$ corresponding to the $x$ component of velocity (right), to the $y$ component of the velocity (middle), and to the pressure (right).}
    \label{fig:StokesVelocity1}
\end{figure}

\begin{table}[bht]
    \centering
    \begin{tabular}{c c c c c }
    & & \multicolumn{2}{c}{Velocity} & \\
    & & ($x$ component) & ($y$ component) & Pressure \\
    \hline
    DD-PGD: & $\frac{\|v^\pgd - v_{\text{ex}}\|_{L^2(\Omega)}}{\|v_{\text{ex}}\|_{L^2(\Omega)}}$ & $8.64 \times 10^{-4}$ & $1.65 \times 10^{-3}$ & $1.32 \times 10^{-3}$\\
    DD-FEM: & $\frac{\|v^h - v_{\text{ex}}\|_{L^2(\Omega)}}{\|v_{\text{ex}}\|_{L^2(\Omega)}}$ & $1.90 \times 10^{-5}$ & $3.81 \times 10^{-5}$ & $1.98 \times 10^{-3}$ \\
    Global FEM: & $\frac{\|v_\Omega^h - v_{\text{ex}}\|_{L^2(\Omega)}}{\|v_{\text{ex}}\|_{L^2(\Omega)}}$ & $4.14 \times 10^{-6}$ & $4.15 \times 10^{-6}$ & $9.87 \times 10^{-4}$ \\
    \hline
    \end{tabular}
    \caption{Relative $L^2(\Omega)$ error for each component of the velocity and the pressure between the exact solution, the DD-PGD solution, the DD-FEM solution and a monolithic finite element solution on the global domain.}
    \label{tab:StokesErrors}
\end{table}
In Figure \ref{fig:StokesVelocity1}, for $\mu =3$, we plot the error $\log_{10} (|v^\pgd(\mu) - v_\text{ex}(\mu)| / \max_\Omega |v_\text{ex}(\mu)|)$, with $v$ equal to the pressure or to either component of the velocity, and we remark that there is no localized loss of accuracy in the overlapping region. (Note that since we compute the $\log_{10}$ of the scaled error, no values are plotted for the velocity components along the top and lateral edges of the domain where the solution is imposed exactly and the nodal error is $0$.)

\begin{figure}[bht]
    \centering
    \includegraphics[width=0.6\textwidth]{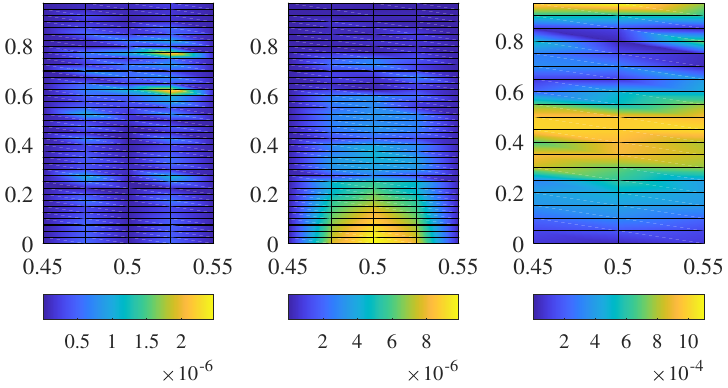}
    \caption{Absolute value of the difference $|v_{\vert\Omega_1}^\pgd - v_{\vert\Omega_2}^\pgd|$ between the local PGD solutions in the overlapping region, for the $x$ component of velocity (left), the $y$ component of velocity (centre), and the pressure (right).}
    \label{fig:Stokes_overlap_diff}
\end{figure}

Finally, in Figure \ref{fig:Stokes_overlap_diff} we present the absolute value in the overlapping region of the difference between the local solutions in $\Omega_1$ and $\Omega_2$ for each of the variables and components. From these plots, it is clear that the coupling algorithm has correctly imposed the continuity of the velocity along the interface (the difference is 0 along the interface edges). Moreover, whilst the pressure coupling has not been directly controlled, there is still good agreement between the local pressure solutions in the overlapping region, with the maximum difference of order $10^{-3}$, which is smaller than the local mesh size. These results agree with those in \cite{Discacciati:2013:Stokes} and, if stronger agreement between the local pressures in the overlapping region is desired, coupling techniques based on the Cauchy stress should be used (see \cite{Discacciati:2013:Stokes}) but, as mentioned in Section \ref{sect:Stokes}, this would increase the computational cost of the online phase.

\subsection{Stokes-Darcy problem with analytic solution in two subdomains}\label{sec:Stokes-Darcy_synthetic}

The second test is a Stokes-Darcy problem with analytic solution that serves two purposes. The first purpose is to demonstrate that the choice of discretization does not impact the quality of the local surrogate models, whilst the second is to assess the performance of the DD-PGD method in the context of a multi-physics problem.

Let $\Omega_1 = [0,1] \times [0.45,1]$ and $\Omega_2 = [0,1] \times [0,0.55]$ be the Stokes and Darcy subdomains, respectively, with Dirichlet and Neumann boundary conditions imposed as shown in Fig.~\ref{fig:Stokes-Darcy_domain-analytical}. We consider a modification of the Stokes-Darcy problem from \cite{Vassilev:2009:SISC}, where the exact velocity and pressure for both problems become
\begin{align}
    \label{eq:SDexact_vel}
    \bu_{\text{ex}}(\bx;\mu_1,\mu_2) &= \begin{bmatrix}
        \phantom{-}\sin\left(\frac{x\mu_1}{\sqrt{\nu K}} + \mu_2\right)e^{\frac{y\mu_1}{\sqrt{\nu K}}}\\
        -\cos\left(\frac{x\mu_1}{\sqrt{\nu K}} + \mu_2\right)e^{\frac{y\mu_1}{\sqrt{\nu K}}}
    \end{bmatrix},\\
    \label{eq:SDexact_pressure}
    p_{\text{ex}}(\bx; \mu_1,\mu_2) &= \sqrt{\frac{\nu}{K}}\left(\frac{1}{\mu_1} - \mu_1\right)\cos\left(\frac{x\mu_1}{\sqrt{\nu K}} + \mu_2\right)e^{\frac{\mu_1}{2\sqrt{\nu K}}} + y - \frac{1}{2}.
\end{align}
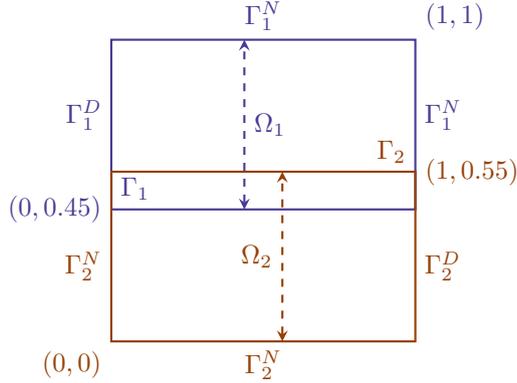
\begin{figure}[bht]
    \centering
    \begin{tikzpicture}
        \draw[BlueViolet, thick] (0, 1.75) rectangle node[pos = 1, above right] {$(1, 1)$} node[pos = 0, left] {$(0, 0.45)$} (4, 4);
        \draw[RawSienna, thick] (0, 0) rectangle node[pos = 0, below left] {$(0, 0)$} node[pos = 1, right] {$(1, 0.55)$} (4, 2.25);
        
        \draw[BlueViolet, stealth-stealth, dashed, thick] (1.75, 1.75) -- node[midway, right]{$\Omega_1$} (1.75, 4);
        \node[BlueViolet, above right] at (0, 1.75) {$\Gamma_1$};
        \node[BlueViolet, left] at (0, 3) {$\Gamma_1^D$};
        \node[BlueViolet, right] at (4, 3) {$\Gamma_1^N$};
        \node[BlueViolet, above] at (2, 4) {$\Gamma_1^N$};

        \draw[RawSienna, stealth-stealth, dashed, thick] (2.25, 0) -- node[midway, left]{$\Omega_2$} (2.25, 2.25);
        \node[RawSienna, above left] at (4, 2.25) {$\Gamma_2$};
        \node[RawSienna, left] at (0, 1) {$\Gamma_2^N$};
        \node[RawSienna, right] at (4, 1) {$\Gamma_2^D$};
        \node[RawSienna, below] at (2, 0) {$\Gamma_2^N$};	
    \end{tikzpicture}
\caption{Computational domain for the Stokes-Darcy problem with analytic solution of Sect. \ref{sec:Stokes-Darcy_synthetic}.}
\label{fig:Stokes-Darcy_domain-analytical}
\end{figure}
The fluid viscosity is $\nu=0.1$, the permeability in the Darcy domain is scalar and equal to $K=1$, while $\mu_1 \in [0.1,1]$ and $\mu_2 \in [1,2]$ are two scalar parameters. The boundary data and forcing terms are reported in \ref{app:StokesDarcyData}. For the space discretization, we consider structured meshes with $h = 2.5 \times 10^{-2}$ both in $\Omega_1$ and $\Omega_2$, and $h_{\mu_1} = 10^{-1}$ and $h_{\mu_2} = 5 \times 10^{-2}$ for the parametric domains of $\mu_1$ and $\mu_2$, respectively. We use stabilized $\mathbb{Q}_1-\mathbb{Q}_1$ elements both for Stokes and Darcy as discussed in Sects.~\ref{sec:StokesProbStatement} and \ref{sec:twoDomainStokesDarcy}. We set $\delta = 5.5$ in the Galerkin Least Squares stabilization for Stokes, while for Darcy, being $K$ a scalar permeability, we slightly modify the definition of the bilinear form $\widetilde{\mathcal{A}}_2$ to include an additional term that accounts for the zero-divergence constraint \cite{Masud:2002:CMAME}. More precisely, we consider
\begin{equation*}
    \widetilde{\mathcal{A}}_2 (\bu,\bv;\bmu) = \frac{1}{2} \int_{\Omega_2} \left( \nu K^{-1} \bu(\bmu) \right) \cdot \bv(\bmu) \, d\bx \\
     + \frac{\beta}{2} \int_{\Omega_2} \nu  K^{-1} \, h^2 \, (\nabla \cdot \bu (\bmu)) \, (\nabla \cdot \bv (\bmu)) \, d\bx\,,
\end{equation*}
with the stabilization coefficient $\beta = 1$.

First, to quantify the impact of the stabilized $\mathbb{Q}_1-\mathbb{Q}_1$ discretization versus the $\mathbb{Q}_2-\mathbb{Q}_1$ used in Sect.~\ref{sect:StokesStokes}, we focus on the Stokes problem only, and we solve the local problems \eqref{eq:StokesLocalPGDData} and the independent problems \eqref{eq:StokesLocalPGDUnity} to construct the PGD expansion \eqref{eq:stokesPGDbasis}. Then, we compute $\bLambda_{\Gamma_1} = \bu_{\text{ex}}(\bx; \mu_1,\mu_2)$ for $\bx \in \Gamma_1$, and we use the values of $\bLambda_{\Gamma_1}$ at the dofs on $\Gamma_1$ to construct the Stokes surrogate models in $\Omega_1$. Considering the Dirichlet boundary condition on $\Gamma_1^D$, there are 40 dofs for each component of the velocity for $\mathbb{Q}_1$ elements and 80 dofs for $\mathbb{Q}_2$ elements. The number of local problems solved in the offline phase for Stokes and the computational cost is reported in Table \ref{table:stokesDiscretisationTest}, which shows that the stabilized $\mathbb{Q}_1-\mathbb{Q}_1$ method leads to a reduction of the computational cost by 2 minutes.

\begin{table}[bht]
\begin{center}
    \begin{tabular}{ccccc}
        Discretization & Problems \eqref{eq:StokesLocalPGDData} & Problems \eqref{eq:StokesLocalPGDUnity} & PGD modes & CPU time  \\
        \hline
        Stabilized $\mathbb{Q}_1-\mathbb{Q}_1$ & 1 & 80 & \phantom{1}94 (139) & 103 s \\
        $\mathbb{Q}_2-\mathbb{Q}_1$ & 1 & 160 & 172 (210) & 223 s \\
        \hline
    \end{tabular}
    \caption{Number of local problems solved, PGD modes and CPU time for the Stokes problem using two space discretizations. (The number of PGD modes in brackets is before compression.)}
    \label{table:stokesDiscretisationTest}
\end{center}
\end{table}

Analogously to Sect. \ref{sect:StokesStokes}, we can compare the accuracy of the two discretization approaches by computing the relative $L^2(\Omega)$ error for each of the variables and components with respect to the known analytical solution for a fixed set of parameter values. The errors computed for $\mu_1 = 0.5$ and $\mu_2 = 1.05$ are presented in Table \ref{tab:StokesDiscretisationComparison}, where we can see that the two methods provide comparable accuracy.

\begin{table}[bht]
    \centering
    \begin{tabular}{c c c c }
        & \multicolumn{2}{c}{Velocity} &  \\
        Discretization & ($x$-component) & ($y$-component) & Pressure\\
        \hline
        Stabilized $\mathbb{Q}_1-\mathbb{Q}_1$& $1.32 \times 10^{-3}$ & $3.35 \times 10^{-3}$ & $9.78 \times 10^{-3}$ \\
        $\mathbb{Q}_2-\mathbb{Q}_1$ & $6.77 \times 10^{-4}$ & $1.79 \times 10^{-3}$ & $1.40 \times 10^{-2}$\\
        \hline
    \end{tabular}
    \caption{Relative $L^2(\Omega)$ error for each component of the velocity and pressure between the exact solution \eqref{eq:SDexact_vel}-\eqref{eq:SDexact_pressure} and the PGD solutions for the Stokes problem taking $\mu_1 = 0.5$ and $\mu_2 = 1.05$.}
    \label{tab:StokesDiscretisationComparison}
\end{table}

Then, we proceed to compute the local surrogate models for Darcy in $\Omega_2$. In this case, in the offline phase, we must solve one problem \eqref{eq:problemPGDDarcyData} and 41 problems like \eqref{eq:problemDarcyBase}. The computational cost is $181$~s to obtain $57$ modes ($113$ modes before compression), that, considering also the Stokes problem, brings the total computational cost of the offline phase for Stokes-Darcy to $284$~s.
The higher computational cost of the offline phase compared to the one for the Stokes-Stokes problem of Sect.~\ref{sect:StokesStokes} is mainly due to the fact that the boundary data and the forces $\beff_i(\bmu)$ that correspond to the exact solution \eqref{eq:SDexact_vel}-\eqref{eq:SDexact_pressure} do not have analytic separated expressions. Hence, the construction of a separated form of the data is performed numerically, e.g., via the \texttt{pgdTensorSeparation} routine described in \ref{app:StokesDarcyData}. In the numerical test of this section, the outcome of the PGD separation is stored to disk, thus requiring loading and reading from disk each time the separated form of the data is evaluated. Whilst this is more costly than the corresponding operations in Sect.~\ref{sect:StokesStokes}, it can be easily optimized by a suitable refactoring of the operations involved in the offline phase, thus retrieving computational times comparable to those of the Stokes-Stokes case.

In the online phase, solving the interface system \eqref{eq:stokesDarcyPGDinterface} for $\mu_1 = 0.5$ and $\mu_2=1.05$ by GMRES requires 33 iterations with a total computational time of 0.21 s ($0.05$ s for only the GMRES iterations). The global Stokes-Darcy surrogate model is then obtained as in \eqref{eq:stokesDarcyPGDglobal}, and the relative errors in $L^2(\Omega)$ norm with respect to the exact solution are $1.37 \times 10^{-3}$, $3.96 \times 10^{-3}$ and $8.45 \times 10^{-3}$ for the $x$-component, $y$-component of the velocity, and for the pressure, respectively. Similarly to Sect.~\ref{sect:StokesStokes}, we notice that the errors agree with the expected finite element accuracy, $O(h^2)$, for the chosen mesh size $h$.

\begin{figure}[ht]
    \centering
    \begin{tabular}{ccc}
      \footnotesize{$x$ component of the velocity}
    & \footnotesize{$y$ component of the velocity}
    & \footnotesize{pressure} \\
    \includegraphics[width=0.3\textwidth]{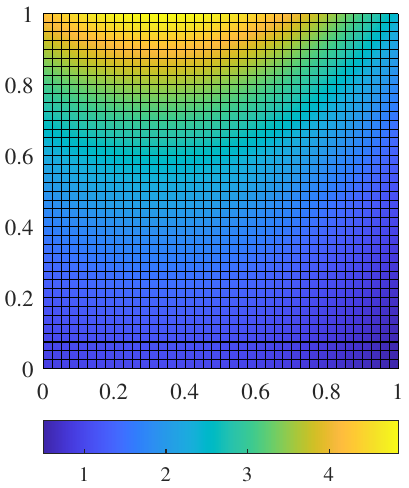}
    & \includegraphics[width=0.3\textwidth]{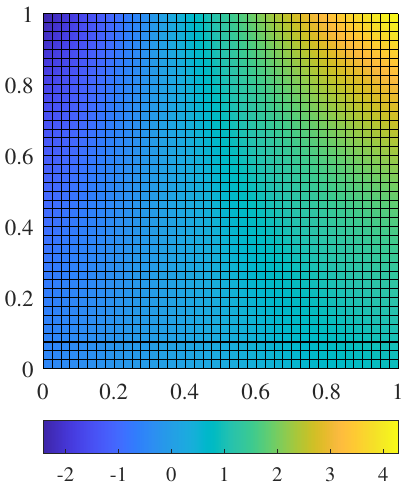}
    & \includegraphics[width=0.3\textwidth]{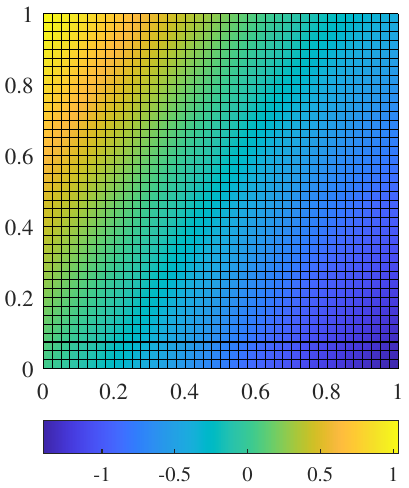}
    \\
    \includegraphics[width=0.3\textwidth]{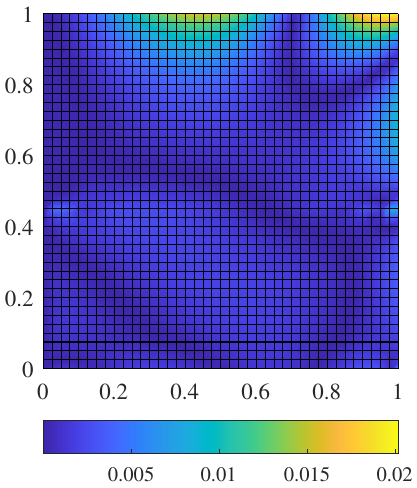}
    & \includegraphics[width=0.295\textwidth]{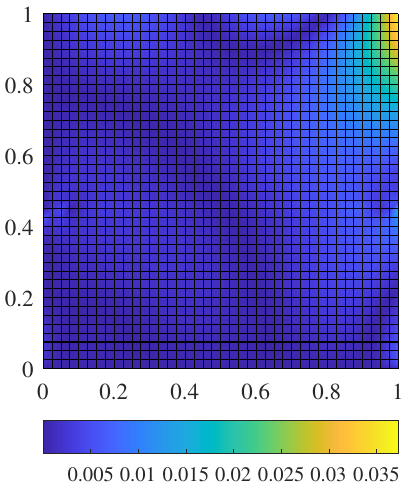}
    & \includegraphics[width=0.295\textwidth]{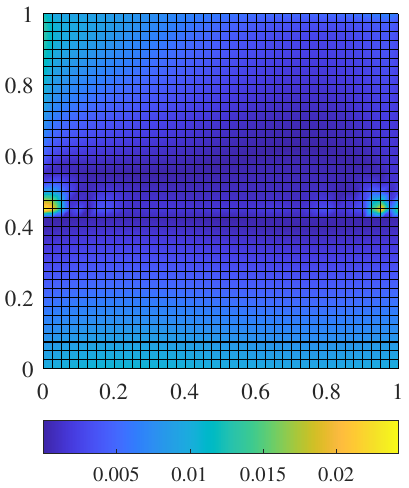} \\
    \end{tabular}
    \caption{The global DD-PGD surrogate model for the parameter configuration $\mu_1 = 0.5$ and $\mu_2 = 1.05$ (top) and the absolute nodal difference between the DD-PGD solution and the analytical solution (bottom) for the $x$-component of the velocity (left), the $y$ component of the velocity (middle), and the pressure (right).}
    \label{fig:stokes-darcy_synthetic}
\end{figure}

In Fig.~\ref{fig:stokes-darcy_synthetic}, we plot the global DD-PGD solution and the absolute value of the pointwise difference between the DD-PGD and the analytical solutions for each variable and component. We remark that, while the errors are all of the same order, the higher values of the error for the pressure are near the interface $\Gamma_1$, where the coupling condition \eqref{eq:couplingConditionSD1} imposes the continuity of the velocity only, so a non-perfect match of the pressures can be expected.

\subsection{Stokes-Darcy membrane cross-flow problem}\label{sect:StokesDarcyPhysical}

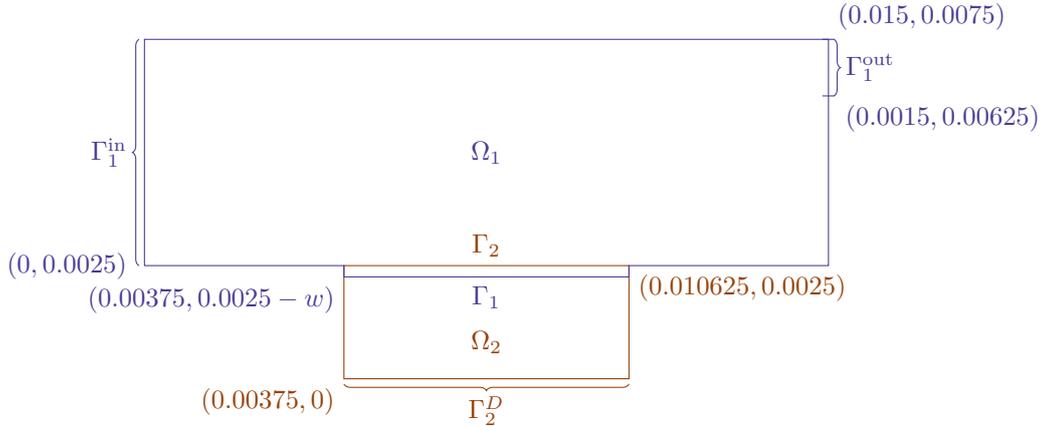
\begin{figure}[bht]
    \centering
    \begin{tikzpicture}
        \draw[RawSienna] (2.625, 0) rectangle node[pos = 0.5, below] {$\Omega_2$} (6.375, 1.5);
        \node[below left, RawSienna] at (2.625, 0) {$(0.00375, 0)$};
        \node[below right, RawSienna] at (6.375, 1.5) {$(0.010625, 0.0025)$};
        \node[above, RawSienna] at (4.5, 1.5) {$\Gamma_2$};
        \draw[RawSienna, decorate, decoration = {brace, raise = 2pt, mirror}] (2.625,0) -- node[below = 3pt] {$\Gamma_2^D$} (6.375, 0);
        
        \draw[BlueViolet] (0, 1.5) -- (2.625, 1.5) -- (2.625,  1.35) -- (6.375, 1.35) -- (6.375, 1.5) -- (9, 1.5) -- (9, 4.5) -- (0, 4.5) -- (0, 1.5);
        \node[BlueViolet] at (4.5, 3) {$\Omega_1$};
        \node[left = 3pt, BlueViolet] at (0, 1.5) {$(0, 0.0025)$};
        \node[above right, BlueViolet] at (9, 4.5) {$(0.015, 0.0075)$};
        \node[below left, BlueViolet] at (2.625, 1.35) {$(0.00375, 0.0025 - w)$};
        \node[below, BlueViolet] at (4.5, 1.35) {$\Gamma_1$};

        \draw[BlueViolet] (8.92, 3.75) -- (9.08, 3.75);
        \node[below right, BlueViolet] at (9.1, 3.75) {$(0.0015, 0.00625)$};
        \draw[BlueViolet] (8.92, 4.5) -- (9.08, 4.5);
        \draw [BlueViolet, decorate, decoration = {brace, mirror, raise = 2pt}] (9, 3.75) -- node[right=3pt] {$\Gamma_1^\text{out}$} (9, 4.5);  

        \draw[BlueViolet, decorate, decoration = {brace, raise = 2pt}] (0, 1.5) -- node[left=3pt] {$\Gamma_1^\text{in}$} (0, 4.5);
    \end{tikzpicture}
    \caption{The physical domain for the Stokes-Darcy crossflow problem.}
    \label{fig:Stokes-Darcy-crossflow-physical-domain}
\end{figure}

In this test case, we simulate the cross-flow of a fluid across a porous medium in two-dimensions, as seen in \cite{Hanspal:2009:ChemEng, Discacciati:2017:IMAJNA}. The domain geometry is illustrated in Fig.~\ref{fig:Stokes-Darcy-crossflow-physical-domain}, with the free flow domain $\Omega_1$ extended by width $w$ into the Darcy domain $\Omega_2$ to introduce an overlap between the subdomains. Gravitational effects are neglected, hence $\beff_i(\bmu)$ is zero in both subdomains. In the Stokes subdomain, $\Omega_1$, we impose the following boundary conditions: on the inlet $\Gamma_1^\text{in} = \{0\} \times (0.0025, 0.0075)$ m, we impose the parametric inflow velocity profile 
\begin{equation}
    \bg_1^D(\mu_1) = 
    (-16000y^2 + 160y - 0.3)
    \begin{bmatrix}
     \cos(\mu_1) \\
     \sin(\mu_1) \\
    \end{bmatrix} \text{ m/s}.
\end{equation}
The parameter $\mu_1 \in [-\pi/2, \pi/2]$ is used to vary the angle of attack of the inflow profile, with the case $\mu_1 = 0$ corresponding to the setting in \cite{Discacciati:2017:IMAJNA}. On $\Gamma_1^\text{out}$, a traction-free condition $\bsigma(\bu(\bmu), p(\bmu); \bmu) \bn = 0$ is imposed to model an outflow boundary. On the remaining boundary $\partial\Omega_1 \setminus (\Gamma_1^\text{in} \cup \Gamma_1^\text{out} \cup \Gamma_1)$, no-slip conditions are enforced.
For the Darcy domain $\Omega_2$, we consider the homogeneous boundary conditions $p_2(\bmu) = 0$ on $\Gamma_2^D = (0.00375, 0.010625) \times \{0\}$ m, and the non-penetration condition $\bu_2(\bmu) \cdot \bn = 0$ on $\partial\Omega_2 \setminus (\Gamma_2^D \cup \Gamma_2)$.

Following \cite{Discacciati:2017:IMAJNA}, we move to a dimensionless setting using the scalings $X_1 = 5 \times 10^{-3}$ m, $U_1 = 10^{-1}$ m/s, and $P_1 = 10$ kg/(m s$^2$) for the space, velocity and pressure, respectively.
Then, we set the dimensionless non-parametric viscosity $\nu = 0.002$, and we introduce the scalar parameter $\mu_2 \in [2 \times 10^{-6},\ 2 \times 10^{2}]$ to represent the scaled inverse scalar parametric permeability $\nu\,K^{-1}$. Remark that, despite the global solution depending on both parameters $\bmu = (\mu_1,\mu_2)$, since we employ domain decomposition, we only need to compute local surrogate models that depend on the parameter that is active in each subdomain, $\mu_1$ in $\Omega_1$ and $\mu_2$ in $\Omega_2$.  
This leads to a reduction in the computational cost of the offline phase, with both Stokes and Darcy problems featuring two spatial and one parametric dimensions.

For the space discretization, the computational mesh is characterized by $h = 1.5625 \times 10^{-2}$, and includes a localized refinement in the vertical direction near the overlapping region in such a way that $w/X_1 = h/4$. Notice that the width the overlapping region is very thin, since this represents the narrow transition region between the free fluid regime and the porous medium regime (see \cite{Discacciati:2024:JCP} for a detailed discussion of the physical interpretation). 
Each component of the velocity and the pressure have 14,178 dofs in $\Omega_1$ and 3,649 dofs in $\Omega_2$, and there are 87 dofs on $\Gamma_1$ and 89 on $\Gamma_2$. Note that the different number of dofs on the interfaces is due to the Dirichlet boundary condition imposed on the edges adjacent to the interface in the Stokes domain. Therefore, for this Stokes-Darcy problem, the number of local problems to solve in the offline phase is 175 in $\Omega_1$ and 89 in $\Omega_2$. The parametric interval for $\mu_1$ is discretized using $h_{\mu_1} = \pi/8$, while the one for $\mu_2$ is replaced by the discrete set of nine points $\mathcal{I}_{\mu_2} = \{ 2 \times 10^i \; : \; i=-6,\ldots,2\}$. In the offline phase, for the Stokes problem we compute a total of 176 modes, with the number of modes being unchanged by PGD compression, meanwhile for Darcy we compute 4,840 modes, reduced to 178 by PGD compression. 
The total computational time is approximately 119 minutes, of which approximately 67 are needed for the Stokes problem.

\begin{figure}
    \centering
    \includegraphics{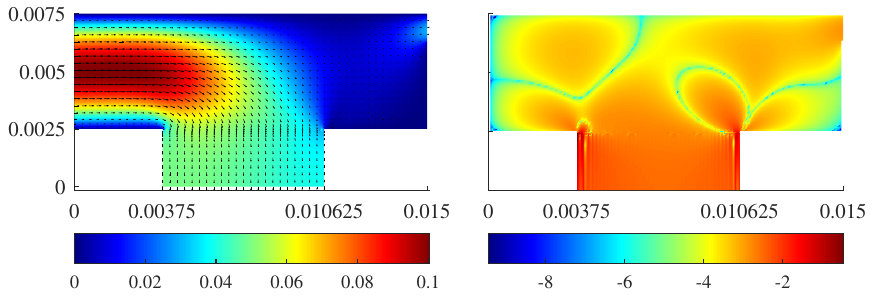}\\
    \includegraphics{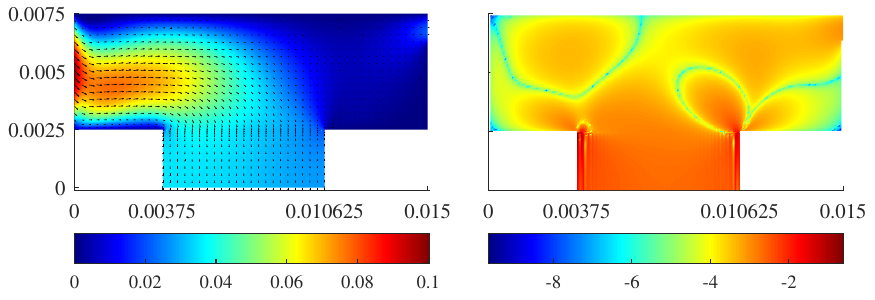}\\
    \includegraphics{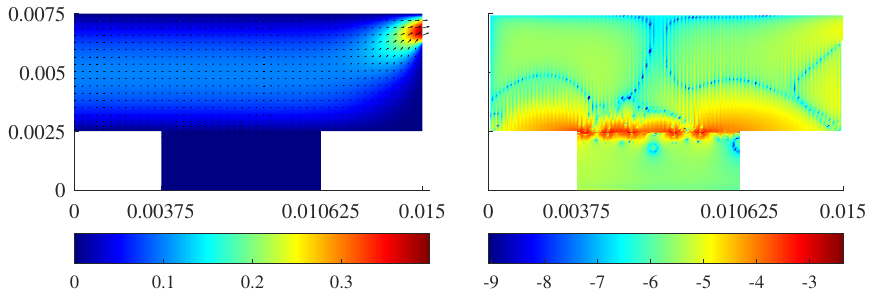}
    \caption{Magnitude of the DD-PGD velocity $\bu^\pgd (\bar{\bmu})$ (left) and error $\log_{10}(| \bu^\pgd(\bar{\bmu}) - \bu^h(\bar{\bmu})|/\max_{\Omega}|\bu^h(\bar{\bmu})|)$ (right) for $\bar{\bmu} = (0, 20)$ (top row), $\bar{\bmu} = (-\pi/4, 20)$ (middle row), and $\bar{\bmu} = (0, 2 \times 10^{-5})$ (bottom row).}
    \label{fig:Stokes-Darcy_crossflow_velocity_1}
\end{figure}

In the online phase, first we consider two parameter configurations that correspond to the case of significant flow into the Darcy domain: $\bar{\bmu} = (0,20)$ and $\bar{\bmu} = (-\pi/4,20)$.
The velocity fields computed for these sets of parameters are shown in the top row and in the middle row of Fig.~\ref{fig:Stokes-Darcy_crossflow_velocity_1}. The computational time of the online phase is only about $0.8$~s for either problem to perform the 144 GMRES iterations needed to solve the interface system \eqref{eq:stokesDarcyPGDinterface}. For this problem due to the high number of iterations the time devoted to the GMRES (approximately $0.75$~s) takes up the bulk of the online phase timing. In comparison, a standard DD-FEM coupling procedure requires 78 iterations and has a computational time of approximately 24 minutes. Thus, the computation of 5 instances of the full-order solver already exceeds the computational cost of the offline phase of the DD-PGD approach.
Then, we consider a final configuration where there is not significant inflow into the porous medium corresponding to $\bar{\bmu} = (0,2 \times 10^{-5})$. This case requires only 7 GMRES iterations for the DD-PGD and 3 for the DD-FEM. But whilst the DD-PGD takes $0.18$~s ($0.08$~s for the GMRES iterations only), DD-FEM requires approximately 2 minutes. The fluid flow simulated for this choice of parameters is plotted in the bottom row of Fig. \ref{fig:Stokes-Darcy_crossflow_velocity_1}.

\begin{figure}
    \centering
    \includegraphics{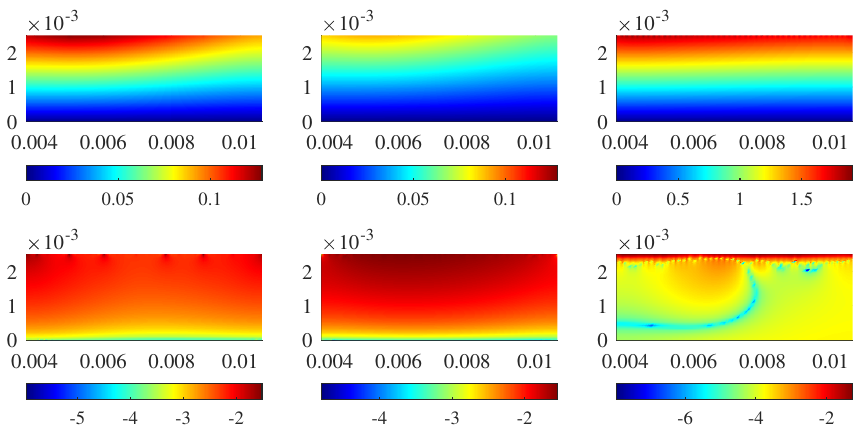}
    \caption{Darcy pressure $p^\pgd(\bar{\bmu})_{\vert\Omega_2}$ (top) and error $\log_{10}(| p^\pgd(\bar{\bmu})_{\vert\Omega_2} - p^h(\bar{\bmu})_{\vert\Omega_2}|/\max_{\Omega_2}|p^h(\bar{\bmu})|)$ (bottom) for $\bar{\bmu} = (0, 20)$ (left column), $\bar{\bmu} = (-\pi/4, 20)$ (right column), and $\bar{\bmu} = (0,2\times10^{-5})$ (right column). 
    }
    \label{fig:Stokes-Darcy_crossflow_pressure_all}
\end{figure}

The corresponding pressure field in $\Omega_2$ for the three choices of the parameters is shown in Fig. \ref{fig:Stokes-Darcy_crossflow_pressure_all}. It is worth noticing that the pressure solutions vary by more than one order of magnitude across the different parameter configurations and yet the local surrogate PGD models can accurately capture this behaviour.

Figures \ref{fig:Stokes-Darcy_crossflow_velocity_1} and \ref{fig:Stokes-Darcy_crossflow_pressure_all} also respectively show the errors $\log_{10}( |\bu^\pgd(\bar{\bmu}) - \bu^h(\bar{\bmu})|/\max_\Omega |\bu^h(\bar{\bmu})|)$ and $\log_{10}(|p^\pgd(\bar{\bmu}) - p^h(\bar{\bmu})|/\max_{\Omega_2} |p^h(\bar{\bmu})|)$, where $\bu^h(\bar{\bmu})$ and $p^h(\bar{\bmu})$ are the finite element solutions computed by the DD-FEM approach on the same mesh used to obtain the local PGD surrogate models. The error plots show that there is very good agreement between the DD-PGD and DD-FEM solutions, with differences of order $10^{-2}$ near the corners of the overlapping region, where the interfaces meet the external boundaries and boundary conditions of different types are imposed.  

\section{Conclusions}
\label{sect:Conclusion}

In this work, a strategy to construct physics-based local surrogate models was presented for parametric Stokes and Stokes-Darcy problems.
The method relies on coupling an overlapping domain decomposition approach to reduce the number of globally coupled degrees of freedom in space with the proper generalized decomposition to reduce the dimensionality of parametric problems.
The DD-PGD method defines arbitrary Dirichlet boundary conditions at the interfaces between the subdomains and constructs low-dimensional local surrogate models exploiting the linearity of the underlying PDEs.
Hence, in the online phase, the subdomains are coupled by solving a parametric interface equation to determine the weighting factors associated with each of the previously computed local surrogate models, without the need to solve any additional PDE in the reduced space.

The method inherits the advantages of the DD-PGD paradigm presented in~\cite{Discacciati:2024:CMAME,Discacciati:2024:DD28}. 
The arbitrary Dirichlet boundary conditions at the interfaces are imposed using the traces of the finite element functions employed for the discretization within the subdomains, avoiding the introduction of new auxiliary basis functions for the solution at the interfaces.
The surrogate models are glued together by enforcing the continuity of the local parametric solutions at the interface, circumventing the imposition of the continuity on the entire overlapping region and the introduction of Lagrange multipliers.

This work goes beyond the state-of-the-art by showing the suitability of the DD-PGD approach to treat physical problems with multiple fields (velocity and pressure), with vector-valued unknowns, and involving multiple physics (viscous incompressible flows and flows in porous media).
The resulting DD-PGD approach is shown to be non-intrusive with respect to the employed solution strategy, exclusively requiring access to the matrices arising from the discretization of the full-order problem, thanks to the framework provided by the Encapsulated PGD Algebraic Toolbox~\cite{Diez:2020:ACME}.
Moreover, the DD-PGD approach is independent of the choice of the full-order discretization. The formulation and the numerical experiments showcase that local surrogate models can be seamlessly constructed using either LBB-compliant finite element pairs or stabilized formulations for incompressible Stokes flows, as well as mixed formulations for the Darcy problem.
Finally, different overlapping DD strategies are employed, from the overlapping Schwarz method of the parametric Stokes-Stokes problem, to the interface control domain decomposition method for the parametric Stokes-Darcy system.

Numerical experiments are presented to assess the accuracy, robustness and computational efficiency of the DD-PGD method for single-physics (Stokes-Stokes) and multi-physics (Stokes-Darcy) problems.
Finally, a cross-flow system commonly employed as benchmark problem for membrane design is employed to showcase the suitability of the method to handle parametric studies in coupled, multi-physics settings of interest in chemical engineering, outperforming standard high-fidelity DD-FEM approaches in terms of computing times.

\paragraph{Acknowledgements} The authors acknowledge funding as follows. MD: EPSRC grant EP/V027603/1. BJE: EPSRC Doctoral Training Partnership grant EP/W523987/1. MG: Spanish Ministry of Science, Innovation and Universities and Spanish State Research Agency \\ MICIU/AEI/10.13039/501100011033 (Grant No. PID2023-149979OB-I00); Generalitat de Catalunya
(Grant No. 2021-SGR-01049); MG is Fellow of the Serra H\'unter Programme of the Generalitat de Catalunya.


\bibliographystyle{elsarticle-num}
\bibliography{references_paper_ddpgd_stokesDarcy}

\appendix

\section{Implementation details}
\label{app:Details}

\subsection{Implementation of multi-variable problems in the Encapsulated PGD Algebraic Toolbox}
\label{app:PGDimplentation}

To clarify how the implementation of the PGD problems \eqref{eq:StokesLocalPGDData}, \eqref{eq:StokesLocalPGDUnity}, \eqref{eq:problemPGDDarcyData} and \eqref{eq:problemDarcyBase} is performed, let us refer to a Stokes model problem with two parameters $\mu_1$ and $\mu_2$, and $\nu(\bmu) = \mu_1 + x \mu_1 \mu_2^2$.

The Neumann data is homogeneous, i.e., $\bg_i^N(\bmu) = \mathbf{0}$, and the source term and the extension of the Dirichlet data are given in separated form as
\begin{align*}
    \beff_i(\bmu) &= \begin{bmatrix}
                        \sin(x) \\
                        \cos(y)
                      \end{bmatrix} 
                    + \begin{bmatrix} 
                        5x\cos(y) \\
                        7
                       \end{bmatrix}\mu_1 \,,\\
    \bg_{\Omega_i}^D(\bmu) &= \begin{bmatrix}
                                5y \\
                                2x
                               \end{bmatrix}\mu_2^2 \,.
\end{align*}
Note that data used for this model problem is chosen for the purposes of clarity of exposition and is not related to any of the numerical tests implemented in Sect.~\ref{sect:numericalResults}.

To compute the local surrogate models, we must construct the matrix 
\begin{equation}
	\label{eq:StokesImpMat}
    \begin{pmatrix}
		\mat{A}_i & \mat{B}_i \\
		\mat{B}_i^T & \mat{0}
	\end{pmatrix}
\end{equation}
where $\mat{A}_i$ and $\mat{B}_i$ are block separated tensors that correspond to the terms in the weak formulation,
\begin{align*}
	(\mat{A}_i)_{m, n} &= \mathcal{A}_i^\pgd(\bvarphi_i^n, \bvarphi_i^m; \bmu) \,, \\
	(\mat{B}_i)_{m, n} &= \mathcal{B}_i^\pgd({\varphi}_i^n, \bvarphi_i^m)\,,
\end{align*}
with $\bvarphi_i^m$ and $\bvarphi_i^n$, and $\varphi_i^n$ the finite element basis functions for velocity and pressure, respectively.
For simplicity of exposition we have assumed that a LBB-stable discretization has been employed. In Table \ref{tab:StokesImplementation}, we identify the finite element matrices and vectors that must be assembled using a finite element library along with the notation we use in Algorithm \ref{alg:StokesImpMat}, which explains how to construct the matrix \eqref{eq:StokesImpMat}. Note that $N^v$ and $N^p$ represent the number of nodes in the finite element discretization of velocity and pressure, respectively. 

\begin{table}[H]
	\begin{tabularx}{\textwidth}{lcX}
		\hline
		\hline\\[-10pt]
		Terms in the   & Algebraic   & Dimensions of algebraic counterpart\\
		weak form      & counterpart & 
		\\[2pt]
		\hline
		\multicolumn{3}{c}{Spatial modes} \\
		\hline\\[-10pt]
		$\displaystyle \int_{\Omega_i} 2 \, \nabla^s \bvarphi^n_i : \nabla^s \bvarphi^m_i \,d\bx$ & $(\mat{A}_i^1)_{m, n}$ & $dN^v \times dN^v$ finite element matrix with constant viscosity $b_\nu^1(\bx) = 1$
		\\[10pt]
		$\displaystyle \int_{\Omega_i} 2x \, \nabla^s \bvarphi^n_i : \nabla^s \bvarphi^m_i \,d\bx$ & $(\mat{A}_i^2)_{m, n}$ & $dN^v \times dN^v$ finite element matrix with space-dependent viscosity  $b_\nu^2(\bx) = x$
		\\[10pt]
		$\displaystyle \int_{\Omega_i} {\varphi}^n_i \nabla \cdot \bvarphi^m_i \,d\bx$ & $(\mat{B}_i^1)_{m, n}$ & $dN^v \times N^p$ finite element matrix
		\\[15pt]
		$\displaystyle 0$ & $(\mat{Z})_{m, n}$ & $N^p \times N^p$ finite element matrix of zeros
		\\[10pt]
        $\displaystyle \int_{\Omega_i} \begin{bmatrix} \sin(x) \\ \cos(y) \end{bmatrix} \cdot \bvarphi_i^n \, d\bx$ & $(\mat{F}_i^1)_n$ & $dN^v \times 1$ finite element vector with space dependent source term  $\bb_{i, f}^1(\bx) = \begin{bmatrix}
            \sin(x) \\ \cos(y)
        \end{bmatrix}$
        \\[10pt]
        $\displaystyle \int_{\Omega_i} \begin{bmatrix} 5x\cos(y) \\ 7 \end{bmatrix} \cdot \bvarphi_i^n \, d\bx$ & $(\mat{F}_i^2)_n$ & $dN^v \times 1$ finite element vector with space dependent source term  $\bb_{i, f}^2(\bx) = \begin{bmatrix}
            5x\cos(y) \\ 7
        \end{bmatrix}$
        \\[10pt]
        $\displaystyle \begin{bmatrix} 5y \\ 2x \end{bmatrix}$ & $\mat{g}^D$ & $dN^v \times 1$ vector of the values of the non-parametric component of the function $\bg_{\Omega_i}^D$ at the nodes used to discretize the velocity
        \\[30pt]
		\hline
		\multicolumn{3}{c}{Parametric modes} \\
		\hline\\[-10pt]
		$\displaystyle 1$ & $\bXi_{\nu, 1}^\text{const}, \bXi_{\nu, 2}^\text{const}$ & $\Nmu \times 1$ vector with discrete values of ones
		\\[10pt]
		$\displaystyle \mu_1$ & $\bXi_{\nu, 1}^\text{linear}$ & $\Nmu \times 1$ vector with discrete values of $\mu_1 \in (\mu_{1,\min},\mu_{1,\max})$
		\\[10pt]
		$\displaystyle \mu_2^2$ & $\bXi_{\nu, 2}^\text{quad}$ & $\Nmu \times 1$ vector with discrete values of $\mu_2^2 \in (\mu_{2,\min}^2,\mu_{2,\max}^2)$
		\\[10pt]
		\hline
		\hline
	\end{tabularx}
	\caption{Differential and algebraic formulation of the terms in the parametric problem for the setup of the Encapsulated PGD Algebraic Toolbox.}
	\label{tab:StokesImplementation}
\end{table}

In Algorithm \ref{alg:StokesImpMat}, we present pseudocode for the construction of the matrix \eqref{eq:StokesImpMat}, following the syntax of the variable type `separated tensor' in the Encapsulated PGD Algebraic Toolbox. Note that concatenation of separated tensors behaves in the same way as standard MATLAB matrix concatenation.

\begin{algorithm}[H]
	\caption{Construction of the matrix \eqref{eq:StokesImpMat} in the Encapsulated PGD Algebraic Toolbox.}\label{alg:StokesImpMat}
	\begin{algorithmic}[1]
		\Require{Spatial matrices ($\mat{A}_i^1, \mat{A}_i^2, \mat{B}_i^1, \mat{Z}$) obtained from a finite element library and parametric vectors ($\bXi_{\nu, 1}^{\text{const}}, \bXi_{\nu, 1}^{\text{linear}}, \bXi_{\nu, 2}^{\text{const}}, \bXi_{\nu, 2}^{\text{quad}}$) constructed via collocation}
        \Statex
		\State{Initialization of separated tensors: \\
			{\small \verb|A_separated = separatedTensor;|} \\
			{\small \verb|B_separated = separatedTensor;|} \\
			{\small \verb|Z_separated = separatedTensor;|} 
		}
        \Statex
		\State{Setup of the spatial modes for $\mat{A}_i$: \\
			{\small \verb|A_separated.sectionalData{1,1} =|} $\mat{A}_i^1${\small \verb|;|} \\
			{\small \verb|A_separated.sectionalData{1,2} =|} $\mat{A}_i^2${\small \verb|;|}
		}
		\State{Setup of the parametric modes for $\mat{A}_i$: \\
			{\small \verb|A_separated.sectionalData{2,1} =|} $\bXi_{\nu, 1}^\text{linear}${\small \verb|;|} \\
			{\small \verb|A_separated.sectionalData{3,1} =|} $\bXi_{\nu, 2}^\text{const}${\small \verb|;|} \\
			{\small \verb|A_separated.sectionalData{2,2} =|} $\bXi_{\nu, 1}^\text{linear}${\small \verb|;|} \\
			{\small \verb|A_separated.sectionalData{3,2} =|} $\bXi_{\nu, 2}^\text{quad}${\small \verb|;|}
		}
		\Statex
		\State{Setup of the spatial modes for $\mat{B}_i$: \\
			{\small \verb|B_separated.sectionalData{1,1} =|} $\mat{B}_i^1${\small \verb|;|}
		}
		\State{Setup of the parametric modes for $\mat{B}_i$: \\
			{\small \verb|B_separated.sectionalData{2,1} =|} $\bXi_{\nu, 1}^\text{const}${\small \verb|;|} \\
			{\small \verb|B_separated.sectionalData{3,1} =|} $\bXi_{\nu, 2}^\text{const}${\small \verb|;|} 
		}
		\Statex
		\State{Setup of the spatial modes for zero matrix: \\
			{\small \verb|Z_separated.sectionalData{1,1} =|} $\mat{Z}${\small \verb|;|}
		}
		\State{Setup of the parametric modes for zero matrix: \\
			{\small \verb|Z_separated.sectionalData{2,1} =|} $\bXi_{\nu, 1}^\text{const}${\small \verb|;|} \\
			{\small \verb|Z_separated.sectionalData{3,1} =|} $\bXi_{\nu, 2}^\text{const}${\small \verb|;|} 
		}
		\Statex
		\State{Creation of the full tensor representing the matrix \eqref{eq:StokesImpMat}: \\
			{\small \verb|LHS = [A_separated, B_separated'; B_separated, Z_separated];|}
		}
        \Statex
		\Ensure{The matrix \eqref{eq:StokesImpMat} in form of a separated tensor \small{\verb|LHS|}}
	\end{algorithmic}
\end{algorithm}

For the right-hand side, we need to construct the vector of separated tensors
\begin{equation}
    \label{eq:StokesImpVector}
    \begin{pmatrix}
        \mat{F}_i \\
        \mat{G}_i
    \end{pmatrix},
\end{equation}
as described in Algorithm \ref{alg:StokesImpVector}.

\begin{algorithm}[H]
	\caption{Construction of the vector \eqref{eq:StokesImpVector} in the Encapsulated PGD Algebraic Toolbox.}\label{alg:StokesImpVector}
	\begin{algorithmic}[1]
		\Require{Spatial matrices and vectors ($\mat{A}_i^1, \mat{A}_i^2, \mat{B}_i^1, \mat{F}_i^1, \mat{F}_i^2, \mat{g}^D$) obtained from a finite element library and parametric vectors ($\bXi_{\nu, 1}^{\text{const}}, \bXi_{\nu, 1}^{\text{linear}}, \bXi_{\nu, 2}^{\text{const}}, \bXi_{\nu, 2}^{\text{quad}}$) constructed via collocation}
        \Statex
		\State{Initialization of separated tensors: \\
			{\small \verb|F_separated = separatedTensor;|} \\
			{\small \verb|G_separated = separatedTensor;|}
		}
        \Statex
		\State{Setup of the spatial modes for $\mat{F}_i$: \\
			{\small \verb|F_separated.sectionalData{1,1} =|} $\mat{F}_i^1${\small \verb|;|} \\
			{\small \verb|F_separated.sectionalData{1,2} =|} $\mat{F}_i^2${\small \verb|;|}\\
            {\small \verb|F_separated.sectionalData{1, 3} = |}
            $-\mat{A}_i^1*\mat{g}^D${\small \verb|;|} \\
            {\small \verb|F_separated.sectionalData{1, 4} = |}
            $-\mat{A}_i^2*\mat{g}^D${\small \verb|;|}
		}
		\State{Setup of the parametric modes for $\mat{F}_i$: \\
			{\small \verb|F_separated.sectionalData{2,1} =|} $\bXi_{\nu, 1}^\text{const}${\small \verb|;|} \\
			{\small \verb|F_separated.sectionalData{3,1} =|} $\bXi_{\nu, 2}^\text{const}${\small \verb|;|} \\
			{\small \verb|F_separated.sectionalData{2,2} =|} $\bXi_{\nu, 1}^\text{linear}${\small \verb|;|} \\
			{\small \verb|F_separated.sectionalData{3,2} =|} $\bXi_{\nu, 2}^\text{const}${\small \verb|;|}\\
            {\small \verb|F_separated.sectionalData{2,3} =|} $\bXi_{\nu, 1}^\text{linear}.* \bXi_{\nu, 1}^\text{const}${\small \verb|;|} \Comment{Elementwise multiplication}\\
            {\small \verb|F_separated.sectionalData{3,3} =|} $\bXi_{\nu, 2}^\text{const}.*\bXi_{\nu, 2}^\text{quad}${\small \verb|;|}\\
            {\small \verb|F_separated.sectionalData{2,4} =|} $\bXi_{\nu, 1}^\text{linear}.*\bXi_{\nu, 1}^\text{const}${\small \verb|;|}\\
            {\small \verb|F_separated.sectionalData{3,4} =|} $\bXi_{\nu, 2}^\text{quad}.*\bXi_{\nu, 2}^\text{quad}${\small \verb|;|}
		}
		\Statex
		\State{Setup of the spatial modes for $\mat{G}_i$: \\
			{\small \verb|G_separated.sectionalData{1,1} =|} $-\mat{B}_i^1*\mat{g}^D${\small \verb|;|}
		}
		\State{Setup of the parametric modes for $\mat{G}_i$: \\
			{\small \verb|G_separated.sectionalData{2,1} =|} $\bXi_{\nu, 1}^\text{const}.*\bXi_{\nu, 1}^\text{const}${\small \verb|;|} \\
			{\small \verb|G_separated.sectionalData{3,1} =|} $\bXi_{\nu, 2}^\text{const}.*\bXi_{\nu, 2}^\text{quad}${\small \verb|;|} 
		}
		\Statex
		\State{Creation of the full tensor representing the vector \eqref{eq:StokesImpVector}: \\
			{\small \verb|RHS = [F_separated; G_separated];|}
		}
        \Statex
		\Ensure{The vector \eqref{eq:StokesImpVector} in form of a separated tensor \small{\verb|RHS|}}
	\end{algorithmic}
\end{algorithm}

The required PGD solution can then be obtained via the call
\begin{equation*}
    u_{i, \pgd}^{0, f}(\bmu) = \verb|pgdLinearSolve(LHS, RHS, <optional settings>)|.
\end{equation*}
For details of the optional settings of \texttt{pgdLinearSolve} we refer to the documentation \cite{Diez:2020:ACME}. However, we would like to briefly discuss the optional setting \texttt{`mask'} as this is essential when working with multi-variable problems.

\subsubsection{Encapsulated PGD Algebraic Toolbox: \texttt{mask} for multi-variable problems}
\label{app:masks}

When using the \texttt{pgdLinearSolve} function provided in the Encapsulated PGD Algebraic Toolbox to solve problems with multiple unknowns, e.g., velocity and pressure, a \texttt{mask} is useful for ensuring that the PGD algorithm converges for all unknowns. A \texttt{mask} is a column vector of boolean values, indicating which part of the spatial solution should be used when considering convergence. One can provide multiple \texttt{masks} by utilizing a matrix structure, and \texttt{pgdLinearSolve} imposes that the convergence criteria must hold for each \texttt{mask} to ensure that, if the magnitudes of the different variables differ significantly, the convergence is guaranteed for each variables.

Suppose that the spatial modes are ordered such that all velocity nodes come first and then the pressure nodes follow, and let $N^v$ and $N^p$ be the number of velocity and pressure nodes respectively. To ensure convergence with respect to both velocity and pressure, the required input \texttt{mask} would be
\begin{equation*}
	\text{\texttt{mask}} = \begin{pmatrix}
		\mat{1}_{N^v} & \mat{0}_{N^v} \\
		\mat{0}_{N^p} & \mat{1}_{N^p}
	\end{pmatrix},
\end{equation*}
where, $\mat{1}_{N^\star}$ and $\mat{0}_{N^\star}$ represent vectors with dimension $N^\star \times 1$ ($\star = v,p$) consisting of 1's and 0's respectively. To use the \texttt{mask}, the call to \texttt{pgdLinearSolve} is modified as
\begin{equation*}
    u_{i, \pgd}^{0, f}(\bmu) = \verb|pgdLinearSolve(LHS, RHS, `mask', mask, <optional settings>)|.
\end{equation*}

\subsection{Data for the Stokes-Darcy problem with analytic solution of Sect.~\ref{sec:Stokes-Darcy_synthetic}}
\label{app:StokesDarcyData}

In this appendix, we indicate the data for the Stokes-Darcy problem with analytic solution studied in Sect.~\ref{sec:Stokes-Darcy_synthetic}. The source term for the Stokes problem is 
\begin{equation*}
    \beff_1(\mu_1, \mu_2) = \begin{bmatrix} \left(\mu_1^2 - 1\right) \sin\left(\frac{\mu_2 \sqrt{\nu K} + x\mu_1}{\sqrt{\nu K}}\right)e^{\frac{\mu_1}{2\sqrt{\nu K}}}\\ 1\end{bmatrix}\,. 
\end{equation*}
The Dirichlet boundary condition on the left edge $\Gamma_1^D$ is imposed through the extension function
\begin{equation*}
    \bg_{\Omega_1}^D(\mu_1, \mu_2) = \begin{bmatrix} \sin(\mu_2)e^{\frac{y\mu_1}{\sqrt{\nu K}}}\\ -\cos(\mu_2)e^{\frac{y\mu_1}{\sqrt{\nu K}}}\end{bmatrix}\,.   
\end{equation*}
The Neumann boundary condition on the top and right edges are respectively calculated to be
\begin{equation*}
    \bg_1^{N, \text{top}}(\mu_1, \mu_2) = \begin{bmatrix}
        \frac{\mu_1\nu\sin\left(\frac{\mu_2 \sqrt{\nu K} + x\mu_1}{\sqrt{\nu K}}\right)e^{\frac{\mu_1}{\sqrt{\nu K}}}}{\sqrt{\nu K}}\\
        \frac{\nu \left(\left(\mu_1^2 - 1\right)e^{\frac{\mu_1}{2\sqrt{\nu K}}} + \mu_1^2 e^{\frac{\mu_1}{\sqrt{\nu K}}} \right)\cos\left(\frac{\mu_2\sqrt{\nu K} + x\mu_1}{\sqrt{\nu K}}\right) + \frac{\mu_1\sqrt{\nu K}}{2}}{\mu_1\sqrt{\nu K}}
    \end{bmatrix}
\end{equation*}
and 
\begin{equation*}
    \bg_1^{N, \text{right}}(\mu_1, \mu_2) = \frac{1}{K}\begin{bmatrix}
        \frac{\nu \left(\left(\mu_1^2 - 1\right)e^{\frac{\mu_1}{2\sqrt{\nu K}}} + \mu_1^2 e^{\frac{y\mu_1}{\sqrt{\nu K}}} \right)\cos\left(\frac{\mu_2\sqrt{\nu K} + \mu_1}{\sqrt{\nu K}}\right) - \left(y - \frac{1}{2}\right)\mu_1\sqrt{\nu K}}{\mu_1\sqrt{\nu K}}\\
        \frac{\mu_1\nu\sin\left(\frac{\mu_2 \sqrt{\nu K} + \mu_1}{\sqrt{\nu K}}\right)e^{\frac{y\mu_1}{\sqrt{\nu K}}}}{\sqrt{\nu K}}
    \end{bmatrix}.
\end{equation*}
The Darcy force term is given by
\begin{align*}
    \beff_2(\mu_1, \mu_2) &= \frac{1}{K}\begin{bmatrix}
        \sin\left(\frac{\mu_2\sqrt{\nu K} + x\mu_1}{\sqrt{\nu K}}\right)\left(\mu_1^2e^{\frac{\mu_1}{2\sqrt{\nu K}}} + \nu e^{\frac{y\mu_1}{\sqrt{\nu K}}} - e^{\frac{\mu_1}{2\sqrt{\nu K}}}\right)\\
        K - \nu\cos\left(\frac{\mu_2\sqrt{\nu K} + x\mu_1}{\sqrt{\nu K}}\right)e^{\frac{y\mu_1}{\sqrt{\nu K}}}
    \end{bmatrix}.
\end{align*}
The Neumann boundary condition $\bu_2(\mu_1, \mu_2) \cdot \bn = g_2^N(\mu_1, \mu_2)$, is enforced through the extension functions
\begin{equation*}
    g_{\Omega_2}^{N, \text{left}}(\mu_1, \mu_2) = \sin(\mu_2)e^{\frac{y\mu_1}{\sqrt{\nu K}}}
\end{equation*}
and
\begin{equation}\label{eq:gnAppendix}
    g_{\Omega_2}^{N, \text{bottom}}(\mu_1, \mu_2) = -\cos\left(\frac{x\mu_1}{\sqrt{\nu K}} + \mu_2\right)\,.
\end{equation}
Since the two edges on which the Neumann boundary condition is imposed are perpendicular to each other, and hence the boundary condition for each edge enforces the velocity for different components, there is no need to construct an extension for each edge that is zero on the other edge. The Dirichlet boundary condition is enforced naturally with
\begin{equation*}
    g_2^D(\mu_1, \mu_2) = \sqrt{\frac{\nu}{K}}\left(\frac{1}{\mu_1} - \mu_1\right)\cos\left(\frac{\mu_1}{\sqrt{\nu K}} + \mu_2\right)e^{\frac{\mu_1}{2\sqrt{\nu K}}} + y - \frac{1}{2}.
\end{equation*}

We remark that the data for this problem is not analytically separable, so we use the script \texttt{pgdTensorSeparation} in the Encapsulated PGD Algebraic Toolbox \cite{Diez:2020:ACME} with tolerance $10^{-4}$ to numerically obtain separated representations needed for the numerical tests.

To this purpose, we must first obtain the 4D grid of spatial and parametric nodes that will be used to define the functions to be separated. This can be done using the MATLAB function \texttt{ngrid} as
\begin{equation*}
    \verb|[X, Y, MU_1, MU_2] = ndgrid(x_nodes, y_nodes, mu_1_nodes, mu_2_nodes)|,
\end{equation*}
where $\verb|x_nodes|$, $\verb|y_nodes|$ are the $x$ and $y$ coordinates of the nodes in the spatial finite element mesh, while $\verb|mu_1_nodes|$ and $\verb|mu_2_nodes|$ are the coordinate of the nodes used to discretize the parametric intervals.

Then, to apply the numerical separation procedure to, e.g., $g_{\Omega_2}^{N,\text{bottom}}$ in \eqref{eq:gnAppendix}, we first define
\begin{equation*}
    \verb|discreteFunction = -cos((X.*MU_1)./(sqrt(nu*K)) + MU_2|,
\end{equation*}
and finally call the function
\begin{equation*}
    \verb|separatedFunction = pgdTensorSeparation(discreteFunction, <optional settings>)|.
\end{equation*}
The output \texttt{separatedFunction} is a separated tensor with 4 rows (2 for the spatial dimensions and 2 for the parameters) in the field \texttt{sectionalData}. The optional settings are similar to those for \texttt{pgdLinearSolve}, and we refer again to the documentation \cite{Diez:2020:ACME} for details.

	
\end{document}